\documentclass[reqno, 11pt]{article}

\usepackage{fullpage}

\usepackage[centertags]{amsmath}
\usepackage{amsmath}
\usepackage{amsfonts}
\usepackage{amssymb}
\usepackage{amsthm}
\usepackage{enumerate}

\usepackage{newlfont}
\usepackage{xcolor}

\usepackage{authblk}

\usepackage{stmaryrd}
\SetSymbolFont{stmry}{bold}{U}{stmry}{m}{n}

\usepackage{bbm}
\include{amsthm_sc}

\usepackage{stmaryrd}
\usepackage{mathrsfs}
\usepackage{dsfont}

\usepackage{fancyhdr}
\usepackage{graphicx}
\usepackage{fancybox}
\usepackage{setspace}
\usepackage{cleveref}

\newtheorem{thm}{Theorem}

\newtheorem{lem}[thm]{Lemma}
\newtheorem{prop}[thm]{Proposition}
\newtheorem{cond}[thm]{Condition}
\theoremstyle{definition}
\newtheorem{defn}[thm]{Definition}

\theoremstyle{remark}
\newtheorem{rem}[thm]{Remark}

\newcommand{\R} {\mathbb{R}}
\newcommand{\C} {\mathbb{C}}

\newcommand{\E} {\mathbb{E}}

\newcommand{\p} {\mathbb{P}}

\DeclareMathOperator{\re}{\mathrm{Re}}
\DeclareMathOperator{\im}{\mathrm{Im}}

\newcommand{\caN}{{\mathcal N}}

\newcommand{\bss}{{\boldsymbol s}}

\newcommand{\bsu}{{\boldsymbol u}}

\newcommand{\bsx}{{\boldsymbol x}}
\newcommand{\bsy}{{\boldsymbol y}}

\newcommand{\bsX}{{\boldsymbol X}}
\newcommand{\bsY}{{\boldsymbol Y}}

\newcommand{\wt}{\widetilde}
\newcommand{\wh}{\widehat}

\newcommand{\beq}{ \begin{equation} }
\newcommand{\eeq}{ \end{equation} }
\newcommand{\beqq}{ \begin{equation*} }
\newcommand{\eeqq}{ \end{equation*} }

\newcommand{\dd}{\mathrm{d}}
\newcommand{\ii}{\mathrm{i}}

\renewcommand{\bss}{\boldsymbol{\sigma}}

\renewcommand{\tt}{\boldsymbol{\tau}}

\newcommand{\hZ}{\hat{Z}}

\newcommand{\red}{d_+}
\newcommand{\led}{d_-}

\newcommand{\n}{n}
\newcommand{\q}{\mathbf{Q}_n}
\newcommand{\qq}{\mathbf{Q}}
\newcommand{\al}{\alpha}
\newcommand{\be}{B}
\newcommand{\aln}{\alpha_{\n}}
\newcommand{\ben}{B_{n}}
\newcommand{\s}{\mathsf{s}}
\newcommand{\h}{\wh H}
\newcommand{\aw}{\wh A}

\newcommand{\dw}{\wh D}
\newcommand{\ew}{\wh E}

\newcommand{\lw}{\wh L}
\newcommand{\et}{\tilde E}

\newcommand{\lt}{\tilde L}
\newcommand{\at}{\tilde A}

\newcommand{\rf}{r_1}
\newcommand{\rs}{r_2}
\newcommand{\mea}{\mu}
\newcommand{\van}{\sigma^2}

\newcommand{\varl}{\mathsf{S}}
\newcommand{\varial}{\mathsf{A}}
\DeclareMathOperator{\TW}{TW}
\DeclareMathOperator{\SSK}{SSK}
\DeclareMathOperator{\BSSK}{BSSK}
\DeclareMathOperator{\MP}{MP}
\DeclareMathOperator{\SC}{SC}
\DeclareMathOperator{\GOE}{GOE}

\numberwithin{equation}{section} 
\numberwithin{thm}{section}

\title{Free energy of bipartite spherical Sherrington--Kirkpatrick model} 

\author{Jinho Baik\footnote{Department of Mathematics, University of Michigan,
Ann Arbor, MI, 48109, USA \newline email: \texttt{baik@umich.edu}}
and Ji Oon Lee\footnote{Department of Mathematical Sciences, KAIST, Daejeon, 34141, Korea
\newline email: \texttt{jioon.lee@kaist.edu}}}

\date{\today}

\begin{document}

\maketitle

\begin{abstract}
We consider the free energy of the bipartite spherical Sherrington--Kirkpatrick model. We find the critical temperature and prove the limiting free energy for all non-critical temperature. We also show that the law of the fluctuation of the free energy converges to the Gaussian distribution when the temperature is above the critical temperature, and to the GOE Tracy--Widom distribution when the temperature is below the critical temperature. The result is universal, and the analysis is applicable to a more general setting including the case where the disorders are non-identically distributed.
\end{abstract}

\section{Introduction} \label{sec:intro}

\subsection{Bipartite SSK}

The Sherrington--Kirkpatrick (SK) model and the spherical Sherrington--Kirkpatrick (SSK) model are disordered systems in which the spin variables are subject to Gibbs probability measures defined by random Hamiltonians. 
They can be thought of as finite-temperature versions of the problem of finding the maximum of a random function on either a hypercube (SK model) or a sphere (SSK model). 
As such, there are significant interests in these models and their generalizations in probability and statistical physics, as well as computer science and social science. 
There is a long history to the subject with many important results. 
We refer to \cite{Guerra2015} and references therein.
 
A natural variation is the case when the spins are divided into two (or more) groups such that the spins of different group interact but the spins of same group do not interact. 
When there are two groups, we are lead to the bipartite system. 

The bipartite spherical Sherrington--Kirkpatrick model (SSK) model is defined as follows. 
Let 
\beq
	S_{n-1}=\{\bsu\in \R^n : \|\bsu\|=\sqrt{n}\}
\eeq
be a sphere in $\R^n$.
Let $N_1$ and $N_2$ be two positive integers and consider two types of spin variables $\bss=(\sigma_1, \cdots, \sigma_{N_1})$ and $\tt=(\tau_1, \cdots, \tau_{N_2})$ on two different spheres, 
\beq
	\bss\in S_{N_1-1}, \qquad \tt\in S_{N_2-1}. 
\eeq
Define the Hamiltonian 
\beq \label{eq:Hn}
	H(\bss, \tt) = \frac1{\sqrt {N}} \sum_{i=1}^{N_1} \sum_{j=1}^{N_2} J_{ij}\sigma_i \tau_j , \qquad N:=N_1+N_2,
\eeq 
where $J_{ij}$ are independent random variables of mean $0$ and variance $1$. 
The bipartite SSK model is defined, for each $\beta>0$, by the Gibbs probability measure
\beq
	P(\bss, \tt)= \frac1{Z_{N_1, N_2}} e^{\beta H(\bss, \tt)}, \qquad (\bss, \tt)\in S_{N_1-1}\times S_{N_2-1}
\eeq
where $\beta$ is called the inverse temperature and $Z_{N_1, N_2}$ is the normalization constant, which is also known as the partition function. 
Note that the probability measure depends on the random variables $J_{ij}$. 

The goal of this paper is to study the free energy $F_{N_1,N_2} (\beta) = N^{-1} \log Z_{N_1,N_2} (\beta)$ as $N_1, N_2\to \infty$. 
For small enough $\beta$, Auffinger and Chen obtained a minimization formula for the limiting free energy in \cite{AuffingerChen2014}. 
We mention that their work applies to more general mixed $(p,q)$-spin Hamiltonians with external fields. 
One of the contributions of this paper is the computation of the limiting free energy for the Hamiltonian~\eqref{eq:Hn} for all $\beta$ other than a critical $\beta_c$. 
We also find the critical inverse temperature $\beta_c$ explicitly. 
When $\beta$ is small, our formula agrees with the result of Auffinger and Chen. 
Another contribution of this paper is the evaluation of the next order term.
We obtain the limiting law of the fluctuations, again for all $\beta\neq\beta_c$. 
We show that the fluctuations are Gaussian for $\beta<\beta_c$, and are given by the Tracy--Widom distribution of random matrix theory for $\beta>\beta_c$. 
In this paper, the disorder parameters $J_{ij}$ are not restricted to Gaussian variables. 

For the usual SK and SSK models, the limiting free energy is given by the Parisi formula \cite{Parisi} and Crisanti--Sommers formula \cite{CrisantiSommers}, which were rigorously proved by Talagrand in \cite{TalagrandParisi, TalagrandParisiSpher}. 
The fluctuations were obtained for $\beta$ below a critical value by first Aizenman, Lebowitz, and Ruelle in \cite{AizenmanLebowitzRuelle} and subsequently in \cite{FrohlochZegarlinski1987, CometsNeveu1995, BovierKurkovaLowes}.
There are several recent results for large $\beta$ and and also for the case with the presence of the external field in \cite{BaikLee, SubagZeitouni2017, ChenDeyPanchenkp2017, ChenSen2017, BaikLee2}. 

\subsection{Multi-species SK}

The bipartite Sherrington--Kirkpatrick (SK) model is defined by the same Hamiltonian~\eqref{eq:Hn} but the spins are now assumed to be on a hypercube,
\beq
	(\bss, \tt)\in  
	\{-1, 1\}^{N_1}\times \{-1, 1\}^{N_2} = \{-1, 1\}^{N_1+N_2}.
\eeq
Note that for the spheres, $S_{N_1-1}\times S_{N_2-1}$ is not equal to $S_{N_1+N_2-1}$. 

The bipartite SK model is a special case of the multi-species Sherrington--Kirkpatrick model. 
The multi-species SK model was introduced in \cite{Barra2015}, and it is defined as follows. 
Let 
\beq \label{eq:HHN}
	H_N^{\text{MS}}(\sigma)= \frac1{\sqrt{N}} \sum_{i,j=1}^N g_{ij}\sigma_i\sigma_j, \qquad \bss=(\sigma_1, \cdots, \sigma_N)\in\{-1,1\}^N.
\eeq
be the usual SK Hamiltonian. 
The disorder parameters $g_{ij}$ are independent centered random variables. 
However, we assume that the variances of $g_{ij}$ are not uniform but they depend on the ``species'' of the index $i$ and $j$. 
Let $\mathcal{S}$ be a finite set independent of $N$ and call the elements of $\mathcal{S}$ species.
Fix a map 
\beq
	s: \{1, \cdots, N \} \to \mathcal{S}.
\eeq	
The value $s(i)$ assigns a species to the index $i$. 
Now we assume that the variance of $g_{ij}$ depends only on the species of $i$ and $j$: 
Let 
\beq 
	\Delta^2=(\Delta^2_{st})_{s,t\in \mathcal{S}}, \qquad \Delta^2_{st}=\Delta^2_{ts}
\eeq
be a symmetric matrix with non-negative entries and we assume that 
\beq \label{eq:gsva}
	\E[ g_{ij}^2] =\Delta_{s(i), s(j)}^2 .
\eeq
Setting $N_s = |\{ i : s(i)=s \}|$, the interesting case is when $\frac{N_s}{N} \to r_s\in (0,1)$ as $N\to\infty$ for each $s\in \mathcal{S}$.

The bipartite SK model is the multi-species SK model when $|\mathcal{S}|=2$ and $\Delta^2= \frac14 \left( \begin{smallmatrix} 0 & 1 \\ 1 & 0 \end{smallmatrix} \right)$.
Note that in this case $\Delta^2$ is not positive-semidefinite.
 
The liming free energy of the multi-species model was studied in \cite{Barra2015} and \cite{Pancheko2015}. 
In \cite{Barra2015}, Barra, Contucci, Mingione, and Tantari obtained a lower bound of the limiting free energy assuming that $\Delta^2$ is positive-semidefinite.
On the other hand, Panchenko obtained an upper bound in \cite{Pancheko2015} for general $\Delta^2$. 
When $\Delta^2$ is positive-semidefinite, the upper bound matches with the lower bound, and hence one obtains the limiting free energy. 
The general case, including the bipartite case, remains an open question; see \cite{Barra2011, Barra2014} for some conjectural formulas for the bipartite SK model.

\subsection{Two multi-species SSK models}

Let us consider a spherical version of multi-species SK model. 
We take the same Hamiltonian as~\eqref{eq:HHN} with same disorder parameters $g_{ij}$ satisfying~\eqref{eq:gsva}. 
Note that if $\bss \in\{-1,1\}^N$, then $\|\bss\|=\sqrt{N}$. 
There are two different natural ways of embedding the hypercube. 
One way is that 
\beq \label{eq:Sone}
	\bss\in S_{N-1}.
\eeq
The other way is that 
\beq \label{eq:Sprodspc}
	\bss \in S_{N_{s_1}-1}\times \cdots \times S_{N_{s_m}-1}
\eeq
where $m$ is the number of species, the set of species is denoted by $\mathcal{S}=\{s_1, \cdots, s_m\}$, and $N_{s_k}$ is the 
number of indices corresponding to the species $s_k$ satisfying $\sum_{k=1}^m N_{s_k}=N$. 
In both cases, $\|\bss\|=\sqrt{N}$. 
Therefore, we have two different multi-species spherical Sherrington--Kirkpatrick models, one with spins on one sphere and the other with spins on a product space of spheres.

The bipartite SSK model we introduced earlier corresponds to a special case of~\eqref{eq:Sprodspc}.
In this paper, we focus only on this model. 
However, using a method similar to this paper, one can study the model with~\eqref{eq:Sone} for bipartite case and also some multi-species cases (possibly not positive-semidefinite $\Delta^2$). 
This ``one-sphere multi-species SSK" model will be considered in a separate paper.

\subsection{Connection to random matrices}

We use a special structure of the Hamiltonian~\eqref{eq:Hn} to study its limiting free energy and the fluctuations.
Setting the matrix $J=(J_{ij})$ and considering $\bss$ and $\tt$ as column vectors, the critical points of the function 
$f(\bss, \tt) =  \bss^T J\tt$ (which is a constant multiple of the Hamiltonian) subject to the constraints $\|\bss\|^2=N_1$ and $\|\tt\|^2=N_2$ satisfy the equations
\beq
	J\tt= \lambda_1 \bss, \qquad J^T\bss=\lambda_2 \tt \,,
\eeq
where $\lambda_1$ and $\lambda_2$ are Lagrange multipliers.
These equations imply that $\bss$ is an eigenvector of the matrix $JJ^T$, $\tt$ is an eigenvector of $J^TJ$, and $\lambda_1\lambda_2$ is an eigenvalue of $J^TJ$ (and also $JJ^T$). 

The matrix $J$ is a random matrix with independent and identically distributed entries. 
The matrix $J^TJ$ is called a (constant multiple of) sample covariance matrix (with null covariance) in statistics and also is said to belong to the Laguerre orthogonal ensemble in random matrix theory \cite{Mehta04, Forrester10, BS2010}.
It is one of the fundamental matrices in random matrix theory. 
The behavior of the eigenvalues of $J^TJ$ (the squared singular values of $J$) in the large dimension limit is well-studied. 
 
There is a more direct connection between the random matrices and the free energy. 
In \cite{KosterlitzThoulessJones, BaikLee}, it was shown that the partition function of the usual SSK model can be expressed as a random single integral. 
In this paper, we obtain a similar result for the bipartite SSK model, but this time the random integral is a double integral; see Lemma~\ref{lem:representation}.
This random double integral involves the eigenvalues of $J^TJ$. 
We analyze the double integral asymptotically using the method of steepest-descent.
The reason that we can apply the method of steepest-descent to the random integral is that even though the eigenvalues are random, their fluctuations about their classical locations are small. 
Precise estimates for the locations of the eigenvalues were obtained recently in random matrix theory. 
The ``rigidity" estimates for the eigenvalues of $J^TJ$ were proved by Pillai and Yin in 2014 \cite{PY} for the case $N_1\neq N_2$. The estimates for the case $N_1=N_2$ follow from \cite{AjankiErdosKruger2014}. 
The rigidity estimates are also crucial in proving the universality in random matrix. They are proved for several other random matrix ensembles  \cite{EYY, EKYY1, BEY}.
Our analysis is applicable to a large class of random double integrals under certain general conditions (including the rigidity condition) on a sequence of random variables. 
We obtain the results for the bipartite SSK model as a special case of a more general asymptotic result for random double integrals. 

The strategy above is an extension to our previous works \cite{BaikLee, BaikLee2} for the SSK model. 
A similar idea was also used in an earlier physics paper \cite{KosterlitzThoulessJones} for a non-rigorous analysis for the limiting free energy. 
An important change from our previous work is that the random integral is a double integral this time.
This change adds significant technical difficulties. 
Even in \cite{BaikLee, BaikLee2}, the asymptotic analysis for large $\beta$ (the low temperature regime) was subtle due to the fact that the critical point in the method of steepest-descent is close to a branch point. 
While we could use a certain symmetry to simplify the situation in the SSK model, we lose the symmetry for the double integral in this paper.
This leads us to a more involved analysis; see Section \ref{sec:low} for more discussions.

\subsection{Organization of the paper}

The paper is organized as follows. 
We state the precise definition of the model and state the main results in Section~\ref{sec:results}. 
The main fluctuation results are Theorem~\ref{thm:main} and Theorem~\ref{thm:eigval}. 
We also state the double integral formula of the partition function. 
The asymptotic analysis of the double integral can be carried out under certain general conditions. 
In Section~\ref{sec:rintegral}, we state these conditions and discuss the critical point for the steepest-descent analysis. 
The asymptotic analysis of the general random double integrals is performed for the high temperature regime in Section~\ref{sec:high} and for the low temperature regime in Section~\ref{sec:low}.
Section~\ref{sec:low} is the most technical part of the paper. 
In Section~\ref{sec:proofeigval}, we prove Theorem~\ref{thm:eigval} using the results of Sections~\ref{sec:high} and~\ref{sec:low}. 
In Section~\ref{sec:ProofFlucFromRMT}, we derive Theorem~\ref{thm:main} from Theorem~\ref{thm:eigval} using results from random matrix theory.
In Section~\ref{sec:non-identical}, we briefly discuss the case where the disorders are non-identically distributed.

\subsubsection*{Acknowledgments}
We would like to thank Antonio Auffinger, Wei-Kuo Chen, Dmitry Panchenko, and Eric Rains for useful conversations. 
The work of Jinho Baik was supported in part by NSF grants DMS-1361782 and DMS-1664692, and the Simons Fellows program.
The work of Ji Oon Lee was supported in part by Samsung Science and Technology Foundation project number SSTF-BA1402-04.

\section{Results} \label{sec:results}

In this section, we define the model precisely and state the results.

\subsection{Definitions} \label{sec:defn}

Let 
\beq
	S_{n-1}=\{\bsu\in \R^n : \|\bsu\|=\sqrt{n}\}
\eeq
be a sphere of radius $\sqrt{n}$ in $\R^n$. 
Let $N_1$ and $N_2$ be positive integers and set
\beq
	N=N_1+N_2. 
\eeq
Let $J=(J_{ij})_{i=1, \cdots, N_1, j=1, \cdots, N_2}$ be an $N_1\times N_2$ matrix with i.i.d. entries of mean $0$ and variance $1$. 
Define the Hamiltonian  
\beq
	H(\bss, \tt) = \frac1{\sqrt N} \sum_{i=1}^{N_1} \sum_{j=1}^{N_2} J_{ij}\sigma_i \tau_j  
	= \frac1{\sqrt N}  \langle \bss, J\tt \rangle, \qquad
	(\bss, \tt) \in  S_{N_1-1}\times S_{N_2-1}.
\eeq 
The free energy of the bipartite SSK model at inverse temperature $\beta$ is defined by 
\beq \label{eq:free energy}
	F_{N_1,N_2} (\beta) = \frac{1}{N} \log Z_{N_1,N_2} (\beta) \,,
\eeq
where the partition function $Z_{N_1,N_2}$ is defined by 
\beq \label{eq:free energy2}
	Z_{N_1,N_2} (\beta) = \int_{S_{N_1-1}} \int_{S_{N_2-1}}  e^{ \beta H(\bss, \tt) } \dd\omega_{N_2}(\tt) \dd\omega_{N_1}(\bss). 
\eeq
Here, $\dd\omega_n(\bsu)$ is the uniform probability measure on the sphere $S_{n-1}$. 

We assume the following for $J$. 
Let $J_{ij}$ be independent random variables such that:
\begin{itemize}
	\item The entries are centered with unit variance, i.e., $\E[J_{ij}] = 0$ and $\E[ J_{ij}^2] = 1$.
	\item For any $i, j$, $\E[J_{ij}^3] = W_3$ and $\E[J_{ij}^4] = W_4$ for some constants $W_3, W_4$.
	\item All moments of $J_{ij}$ are finite.
\end{itemize}

We consider the limit as $N, N_1, N_2\to \infty$.
Assume that there is $\delta>0$ such that 
\beq \label{eq:def_rf_rs}
	\frac{N_1}{N} = \rf + O(N^{-1-\delta}), \qquad \frac{N_2}{N}= \rs+ O(N^{-1-\delta}) .
\eeq
for some 
\beq
	\rf, \rs > 0 \,\text{ satisfying }\, r_1 +r_2 = 1.
\eeq

\subsection{Limiting free energy}

We first state the limiting free energy.

\begin{thm} \label{thm:LLN}
Set 
\beq
	\beta_c:= (\rf\rs)^{-\frac14}.
\eeq
Define, for $0<\beta<\beta_c$, 
\beq
	F(\beta)= \frac{\rf\rs\beta^2}{2}
\eeq
and for $\beta>\beta_c$,
\beq \begin{split}
	F(\beta) &=
	\frac{(\sqrt{\rf}+\sqrt{\rs})\sqrt{\varl}-\sqrt{\rf\rs}-1}{2}  \\
	&\qquad \qquad - \frac{\rf-\rs}4 \log \left( \frac{\sqrt{\varl}+\sqrt{\rf}-\sqrt{\rs}}{\sqrt{\varl}-\sqrt{\rf}+\sqrt{\rs}} \right) - \frac{\rs}4 \log \rf - \frac{\rf}{4} \log \rs - \frac12 \log \beta
\end{split} \eeq
where
\beq \label{eq:varlde}
	\varl=\varl(\beta, \rf, \rs):= (\sqrt{\rf}-\sqrt{\rs})^2+4 \rf\rs\beta^2.
\eeq
Then,
\beq
	F_{N_1, N_2}(\beta) \to F(\beta)
\eeq
as $N \to \infty$ in probability for every $\beta\neq \beta_c$. 
\end{thm}

\begin{proof}
This result is a simple consequence of Theorem~\ref{thm:main} below for the fluctuations. 
\end{proof}

Auffinger and Chen obtained the limiting free energy when $\beta$ is small enough in \cite{AuffingerChen2014} in terms of a minimization problem. 
Their result applies to general mixed $(p,q)$-spin Hamiltonians with the presence of the external field. 
The specialization to the $(p,q)=(1,1)$ case (we also set $h_1=h_2=0$ and $\beta_{1,1}=\sqrt{\rf\rs}\beta$ in Theorem 1 of \cite{AuffingerChen2014}) is the following: There is a small constant $\beta_0>0$ (which is not explicitly determined) such that for $\beta<\beta_0$,  
\beq \label{eq:ACformula}
	\lim_{N\to \infty} F_{N_1, N_2}(\beta)= \min_{a,b\in [0,1)}  P(a,b)
\eeq	
where  
\beq
	P(a,b) = \frac{\rf}{2} \left( \frac{a}{1-a} + \log (1-a) \right) + \frac{\rs}{2} \left( \frac{b}{1-b} + \log (1-b) \right)  
	+ \frac{\rf\rs \beta^2}{2} \left( 1-ab \right) .
\eeq
It is easy to find the minimum explicitly. 	
It is straightforward to check that the minimum occurs on the boundary of domain $[0,1)\times [0,1)$ when $\beta\le (\rf\rs)^{-1/4}$ and inside the domain $[0,1)\times [0,1)$ when $\beta>(\rf\rs)^{-1/4}$.
The minimizers are $(a,b)=(0,0)$ when $\beta\le (\rf\rs)^{-1/4}$ and 
\beq
	(a,b)= \left( 1- \frac{\sqrt{\varl} - \sqrt{\rf}+\sqrt{\rs}}{2\sqrt{\rf} \rs \beta^2}, 1- \frac{\sqrt{\varl} + \sqrt{\rf} - \sqrt{\rs}}{2\rf\sqrt{\rs} \beta^2} \right) 
\eeq	
when $\beta>(\rf\rs)^{-1/4}$, where $\varl$ is~\eqref{eq:varlde}. From this, we find that the minimum is equal to $F(\beta)$ 
in Theorem~\ref{thm:LLN} for all $\beta\neq \beta_c$. 
Hence, Theorem~\ref{thm:LLN} implies that the result~\eqref{eq:ACformula} of Auffinger and Chen actually holds for all $\beta\neq \beta_c$ for the $(1,1)$-spin Hamiltonian.

\subsection{Fluctuations of the free energy}

Next result is about the fluctuations of the free energy. 

\begin{thm} \label{thm:main}
We have the following convergence in distribution. 
\begin{enumerate}[(i)]
\item In the high temperature regime $0 < \beta < (\rf\rs)^{-\frac{1}{4}}$, 
\beq
	N(F_{N_1, N_2} - F(\beta)) \Rightarrow \caN(\mea, \van),
\eeq
where $\caN(\mea, \van)$ is the Gaussian distribution with mean
\beq
	\mea = \frac14 \log (1-\rf\rs\beta^4) - \log 2 - (W_4-3) \frac{\rf\rs\beta^4}{4} 
\eeq
and variance 
\beq
	\van = - \frac12 \log \left( 1-\rf\rs\beta^4 \right) + (W_4-3) \frac{\rf\rs\beta^4}{4} .
\eeq

\item In the low temperature regime $\beta >(\rf\rs)^{-\frac{1}{4}}$, 
\beq
	\frac{N^{\frac23}}{\varial}(F_{N_1, N_2} - F(\beta)) \Rightarrow \TW
\eeq
where
\beq \begin{split}
	\varial = \varial(\beta, \rf, \rs)
	&= \frac{ (\sqrt{\rf}+\sqrt{\rs})^{\frac13} ( \sqrt{\varl} - \sqrt{\rf}-\sqrt{\rs}) }{ 4(\rf\rs)^{\frac16} } 
\end{split} \eeq
with $\varl=(\sqrt{\rf}-\sqrt{\rs})^2+4 \rf\rs\beta^2$ defined in~\eqref{eq:varlde} and $\TW$ denotes the GOE Tracy--Widom distribution.
\end{enumerate}
\end{thm}

We remark that the limiting free energy $F(\beta)$ and the constants $\mu$, $\sigma^2$, and $\varial$ are all symmetric in $\rf$ and $\rs$.

The above change from the Gaussian distribution for high temperature to the Tracy--Widom distribution for low temperature also occurs in the usual SSK model  \cite{BaikLee}. 

\subsection{Free energy and eigenvalues}

Assume, without loss of generality, that 
\beq
	N_1\ge N_2.
\eeq 
The matrix of the disorder parameters $J=(J_{ij})$ is an $N_1\times N_2$ matrix. 
We consider the $N_2\times N_2$ square random matrix 
\beq \label{eq:Qmatrx}
	S= \frac1{N_1} J^T J .
\eeq
In statistics, $S$ is known as a sample covariance matrix (with null covariance). 
In random matrix theory, $S$ is also known to belong to the Laguerre orthogonal ensemble \cite{Mehta04, Forrester10, BS2010}. 
Let 
\beq
	\mu_1\ge \mu_2 \ge \dots\ge \mu_{N_2}\ge 0
\eeq
be the eigenvalues of $S$.
We note that $\sqrt{\mu_i}$ are the singular values of $\frac1{\sqrt{N_1}}J$. 

The eigenvalues of $S$ are well studied in the random matrix theory.
For example, the empirical spectral distribution (ESD) of $S$ converges to the Marchenko-Pastur distribution \cite{Marcenko-Pastur63}: 
\beq \label{eq:ESDMP}
	\frac1{N_2} \sum_{i=1}^{N_2} \delta_{\mu_i}(x) \dd x \to \dd\mu_{\MP}(x)
\eeq
weakly in probability as $N_1, N_2\to \infty$ with $\frac{N_2}{N_1} \to \frac{\rf}{\rs} \in (0, 1]$, where
\beq \label{eq:MPformula}
	\dd\mu_{\MP}(x) := \frac{2 \sqrt{(\red-x)(x-\led)}}{\pi (\sqrt{\red}-\sqrt{\led})^2 x} \mathds{1}_{(\led, \red)}(x) \, \dd x
\eeq
with
\beq
	\led= \frac{(\sqrt{\rf}-\sqrt{\rs})^2}{\rf}, \qquad \red= \frac{(\sqrt{\rf}+ \sqrt{\rs})^2}{\rf} .
\eeq

The next theorem relates the second leading term of the free energy with the eigenvalues of $S$. We begin by introducing a suitable notion for the estimates.
\begin{defn}[High probability event]
We say that an $N$-dependent event $\Omega_N$ holds with high probability if, for any given $D > 0$, there exists $N_0 > 0$ such that
$$
\p (\Omega_N^c) \leq N^{-D}
$$
for all $N > N_0$.
\end{defn}

\begin{thm} \label{thm:eigval}
Without loss of generality, assume that $\rf\ge \rs$. 
The following hold with high probability for any fixed $0<\epsilon<\frac{1}{100}$.
\begin{enumerate}[(i)]
\item In the high temperature regime $0 < \beta < (\rf\rs)^{-\frac14}$,  
\beq \begin{split} \label{eq:eigvalh}
	F_{N_1,N_2}(\beta)  
		= F(\beta) - \frac{1}{2N}& \left[ \sum_{i=1}^{N_2} \log  ( z_c - \mu_i) -  N_2 \int \log ( z_c-x) \dd\mu_{\MP}(x) \right]  \\
		+ \frac1{N}& \left[  \frac12 \log \left(1-\rf\rs\beta^4\right) -\log 2 \right]
		+ O(N^{-2+\epsilon}) 
\end{split} \eeq
where 
\beq
	z_c	= \frac{1+\beta^2+\rf\rs\beta^4}{\rf\beta^2}.
\eeq

\item In the low temperature regime $\beta > (\rf\rs)^{-\frac14}$,
\beq \begin{split} \label{eq:Ffllo}
	F_{N_1,N_2}(\beta) &= F(\beta) +   
	\left(\mu_1- \frac{(\sqrt{\rf}+\sqrt{\rs})^2}{\rf} \right) \frac{\rf ( \sqrt{\varl} - \sqrt{\rf}-\sqrt{\rs})}{4(\sqrt{\rf}+\sqrt{\rs})} + O(N^{-1+\epsilon}).
\end{split} \eeq
with high probability where $\varl=(\sqrt{\rf}-\sqrt{\rs})^2+4 \rf\rs\beta^2$ as in~\eqref{eq:varlde}.
\end{enumerate}
\end{thm}

Theorem~\ref{thm:eigval} shows that the difference $F_{N_1,N_2}(\beta) - F(\beta)$ is governed by the top eigenvalue $\mu_1$ when $\beta>\beta_c$ and by a certain combination of all eigenvalues when $\beta<\beta_c$. 
The behaviors of the top eigenvalue and the special combination of all the eigenvalues appearing in the theorem are well-known in random matrix theory. 
In Section~\ref{sec:ProofFlucFromRMT}, we prove Theorem~\ref{thm:main} by combining Theorem~\ref{thm:eigval} and the results from random matrix theory.

\subsection{Special case}

When $\rf=\rs=\frac12$, the formulas are particularly simple. We will compare the formulas with the usual SSK model:
\beq
\begin{split}
	F^{\SSK}(\beta)
	=  \begin{cases} \beta^2 \qquad &\text{for $\beta<\frac12$,} \\
	 2\beta -\frac34   - \frac{1}2 \log \left( 2\beta \right) \quad &\text{for $\beta>\frac12$.}
	 \end{cases}
\end{split}
\eeq 
When $\rf=\rs=\frac12$, we find that the limiting free energy of the bipartite SSK models satisfies 
$$F^{\BSSK}(2\sqrt{2}\beta)= F^{\SSK}(\beta), \qquad \beta\neq\frac12.$$
For general $\rf\neq \rs$, we have 
$$ F^{\BSSK}\left( \beta \sqrt{\frac{2}{\rf\rs}}  \right) = F^{\SSK}(\beta), \qquad \beta<\frac12,$$
but this relationship is not true in low temperature regime $\beta>\frac12$. 

For the fluctuations, when $\rf=\rs$ and $\beta>\frac12$, $\varial(2\sqrt{2}\beta)  =  \beta -\frac12$.
This is the same constant appearing for the low temperature fluctuations of the usual SSK model \cite{BaikLee} (see (iv) of Section 3.1). 
However, when $\rf=\rs$ and $\beta<\frac12$, the constants $\mea (2\sqrt{2}\beta)$ and $\van (2\sqrt{2}\beta)$ are not same as the constants for the high temperature fluctuations of the usual SSK model ((3.12) and (3.13) of \cite{BaikLee}).

We note that the limiting distribution of the eigenvalues associated to the bipartite SSK and the usual SSK are related when $\rf=\rs$. 
When $\rf=\rs$, then the Marchenko-Pastur distribution~\eqref{eq:MPformula} is
\beq
	\mu_{\MP}(x) = \frac{\sqrt{x(4-x)}}{2\pi x} \mathds{1}_{(0, 4)}(x) \, \dd x \,.
\eeq
After a simple change of variables $x = y^2$, this distribution is equal to the semicircle distribution,
\beq
	\mu_{\SC}(y) = \frac{\sqrt{4- y^2}}{2\pi} \mathds{1}_{(-2, 2)}(y) \, \dd y \,,
\eeq
which is the limiting distribution for the random symmetric matrix associated to the usual SSK model.

\subsection{Double integral representation}

As mentioned in Introduction, the starting point of our analysis for Theorem~\ref{thm:eigval} is an explicit double integral formula for the partition function. 
In this subsection, we state and prove the formula. 
Recall that we assume, without loss of generality, that $N_1\ge N_2$. Let 
\beq
	S^{n-1}=\{\bsu\in \R^n : \|\bsu\|=1\}
\eeq 
and $\Omega_n(\bsu)$ is the surface measure (which is not normalized) on the unit sphere $S^{n-1}$. 
After setting $\bss=\sqrt{N_1}\bsx$ and $\tt=\sqrt{N_2}\bsy$, the partition function~\eqref{eq:free energy2} satisfies 
\beq \label{eq:partition}
	Z_{N_1,N_2}(\beta) = \frac{\hat{Z}_{N_1,N_2}(N_1N_2^{\frac12}N^{-\frac12} \beta)}{|S^{N_1-1}||S^{N_2-1}|}
\eeq
where
\beq
	\hat{Z}_{N_1,N_2} (b)= \int_{S^{N_1-1}} \int_{S^{N_2-1}} e^{b \langle \bsx,  M  \bsy \rangle} \dd\Omega_{N_2}(\bsy)\, \dd\Omega_{N_1}(\bsx) , 
	\qquad M:= \frac{J}{\sqrt{N_1}}.
\eeq
Let $\mu_1\ge \mu_2 \ge \dots\ge \mu_{N_2}\ge 0$ be the eigenvalues of the $N_2\times N_2$ matrix $S = M^T M
= \frac1{N_1} J^TJ$. 
The following formula is a variation of a result in \cite{BaikLee}.
 
\begin{lem} \label{lem:representation}
For $N_1\ge N_2$, we have
\beq
	\hZ_{N_1,N_2}(b)= - 2^{N_2} \left( \frac{\pi}{b} \right)^{\frac{N_1+N_2}2-2}
		\int_{\gamma_1-\ii\infty}^{\gamma_1+\ii\infty} \int_{\gamma_2-\ii\infty}^{\gamma_2+\ii\infty} 
	\frac{e^{b (z_1+ z_2)}}{z_1^{(N_1-N_2)/2}\prod_{i=1}^{N_2} \sqrt{4z_1z_2-\mu_i}}  \, \dd z_2 \, \dd z_1
\eeq
where $\gamma_1$ and $\gamma_2$ are any real positive constants satisfying $4\gamma_1\gamma_2>\mu_1$. 
\end{lem}

\begin{proof}
From the singular value decomposition, $M=UD V$ where $U$ and $V$ are orthogonal matrices (of size $N_1$ and $N_2$, respectively) and $D=(D_{ij})$ is an $N_1\times N_2$ matrix with $D_{ii}=\sqrt{\mu_i}$ and $D_{ij}=0$ for $i\neq j$. 
Hence, after changing the variables $\bsx$ and $\bsy$ to $U\bsx$ and $V^T\bsy$, respectively, we have 
\beq
	\hZ_{N_1,N_2} (b)= \int_{S^{N_1-1}} \int_{S^{N_2-1}} e^{b \sum_{i=1}^{N_2} \sqrt{\mu_i} x_i y_i} \, \dd\Omega_{N_2}(\bsy) \dd\Omega_{N_1}(\bsx).
\eeq
Consier 
\beq
	I(z_1, z_2) := \int_{\R^{N_1}} \int_{\R^{N_2}} e^{\sum_{i=1}^{N_2} \sqrt{\mu_i} X_i Y_i} e^{-z_1 |\bsX|^2-z_2 |\bsY|^2} 
	\, \dd^N \bsY  \, \dd^M \bsX .
\eeq
We evaluate this integral in two ways. First, by computing the Gaussian integrals, 
\beq
\begin{split}
	I(z_1, z_2) 
	&= \frac{2^{N_2}\pi^{\frac{N_1+N_2}2}}{z_1^{\frac{N_1-N_2}2}} \prod_{i=1}^{N_2} \frac{1}{\sqrt{4z_1z_2-\mu_i}}
\end{split}
\eeq
for $z_1$ and $z_2$ satisfying $\re z_1>0$, $\re z_2>0$ and $\re(4z_1z_2)>\mu_1$. 
Second, using polar coordinates $\bsX=\sqrt{u}\bsx$, $Y=\sqrt{v}\bsy$ with $u, v>0$ and $\bsx\in S^{N_1-1}$, $\bsy\in S^{N_2-1}$, we find that 
\beq
	I(z_1, z_2) = \int_0^\infty \int_0^\infty \frac14 u^{\frac{N_1}{2}-1}v^{\frac{N_2}{2}-1} \hZ_{N_1,N_2}(\sqrt{uv}) e^{-z_1u-z_2v} \, \dd v \, \dd u .
\eeq	 
By taking the inverse Laplace transform twice, we find 
\beq
	u^{\frac{N_1}2-1} v^{\frac{N_2}2-1} \hZ_{N_1,N_2}(\sqrt{uv}) 
	= -\frac1{\pi^2} \int_{\gamma_1-\ii\infty}^{\gamma_1+\ii\infty} \int_{\gamma_2-\ii\infty}^{\gamma_2+\ii\infty} 
	I(z_1, z_2) e^{z_1u+z_2v} \dd z_2 \, \dd z_1
\eeq
for any $\gamma_1, \gamma_2>0$ satisfying $4\gamma_1\gamma_2>\mu_1$. Setting $u=v=b$, we obtain the result. 
\end{proof}

\section{Random double integral} \label{sec:rintegral}

The main technical part of this paper is the asymptotic analysis of the double integral in Lemma~\ref{lem:representation}.
The integrand contains the random eigenvalues $\mu_i=\mu_i(N_2)$, $1\le i\le N_2$.
We use the method of steepest-descent to evaluate the double integral asymptotically. 
This is possible since the eigenvalues satisfy certain rigidity estimates \cite{PY, AjankiErdosKruger2014} with high probability. 
Since the analysis depends only on the rigidity estimates and a few other properties of $\mu_i$, we present the analysis for a general sequence of random double integrals.  
In this section, we define general random double integrals and state the conditions for the parameters and random variables of the integrals. 
The asymptotic analysis is carried out in the next two sections, Sections~\ref{sec:high} and~\ref{sec:low}. 
Section~\ref{sec:low} is the most technical part of the analysis. 
We then discuss in Section~\ref{sec:proofeigval} that the eigenvalues of the matrix $\frac1{N_1} J^TJ$ for the bipartite SSK model satisfy the conditions (with high probability) and derive Theorem~\ref{thm:eigval} from the general asymptotic results, Proposition~\ref{prop:integral high} and~\ref{prop:lowresult} for the double integrals.

\subsection{General conditions for random double integrals}

Let us define a sequence of general random double integrals.

\begin{defn} \label{defn:integral}
Suppose that for each positive integer $\n$, there are $\n$ non-negative numbers 
$\mu_1(\n)\ge \cdots \ge \mu_{\n}(\n)\ge 0$. 
Let $\aln\ge 0$ and $\ben>0$ be real numbers. 
For each positive integer $\n$, define
\beq \label{eq:integralq}
	\q=\qq(\n,\aln, \ben) := - 
	\int_{\gamma_1-\ii\infty}^{\gamma_1+\ii\infty} \int_{\gamma_2-\ii\infty}^{\gamma_2+\ii\infty} 
	\frac{e^{\n \ben (z_1+ z_2)}}{z_1^{\n \aln}\prod_{i=1}^{\n} \sqrt{4z_1z_2-\mu_i(n)}}  \, \dd z_2 \, \dd z_1
\eeq
where $\gamma_1$ and $\gamma_2$ are any real positive constants satisfying $4\gamma_1\gamma_2>\mu_1(\n)$. 
\end{defn}

We consider large $\n$ asymptotics of $\q$ under the following three conditions.

\begin{cond} \label{cond:aanndb}
There is $\delta>0$ such that 
\beq
	\aln= \al + O(n^{-1-\delta}), \qquad \ben= \be+O(n^{-1-\delta})
\eeq
for $\al\ge 0$ and $\be>0$. 
\end{cond}

\begin{cond}[Regularity of measure] \label{cond:empirical}
The empirical spectral distribution converges weakly in probability to a probability measure $\wh \mu$, i.e.,
\beq \label{eq:ESD_regular}
	\frac1{\n} \sum_{i=1}^{\n} \delta_{\mu_i(\n)}(x) \dd x \to \dd \wh\mu(x),
\eeq
and $\wh\mu$ satisfies the following properties:
\begin{itemize}
 \item $\wh\mu$ is supported on a closed interval $[\led, \red]$.
 \item $\wh\mu$ has a density that is positive on $(\led, \red)$.
 \item The density of $\wh\mu$ exhibits square-root decay at the upper edge, i.e., for some $c_{\wh\mu}>0$, 
\end{itemize}
\beq
	\frac{\dd \wh\mu}{\dd x}(x) = c_{\wh\mu} \sqrt{\red -x} (1+ o(1)) \qquad \text{as } x \uparrow \red.
\eeq
\end{cond}

\begin{cond}[Rigidity] \label{cond:rigidity}
For a positive integer $k \in [1, n]$, let $\hat k := \min \{ k, n +1-k \}$.
Let $g_k$ denote the ``classical location'' defined by the quantiles, 
\beq \label{eq:classicallocationdef}
	\int_{g_k}^{\infty} \dd \wh\mu = \frac{1}{n} \left( k - \frac{1}{2} \right).
\eeq
Then, for any $\epsilon > 0$, 
\beq \label{rigidity}
	| \mu_k(\n) - g_k | \leq \hat k^{-1/3} \n^{-2/3 + \epsilon}
\eeq
for all $1\le k\le \n$ and for all $\n$.
\end{cond}

Note that the last two conditions imply that 
\beq \label{eq:mu1tored}
	\mu_1 (\n) \to \red. 
\eeq

\begin{rem}[Notational Remark 1]
Throughout the paper we use $C$ or $c$ in order to denote a constant that is independent of $\n$. 
Even if the constant is different from one place to another, we may use the same notation $C$ or $c$ as long as it does not depend on $\n$ for the convenience of the presentation. \end{rem}

\begin{rem}[Notational Remark 2]
We use standard notations $O(\cdot)$, $o(\cdot)$, $\ll$, and $\gg$ as $\n \to \infty$. 
\end{rem}

In terms of the above notation, the partition function is given by (see Lemma~\ref{lem:representation}) 
\beq \label{eq:partition_Q}
	Z_{N_1,N_2}(\beta) 
	=  \frac{\qq \big( N_2, \frac{N_1-N_2}{2N_2}, \frac{N_1}{\sqrt{N_2N}} \beta \big)}{|S^{N_1-1}||S^{N_2-1}|} 2^{N_2} \left( \frac{\pi^2 N}{N_1^2N_2 \beta^2} \right)^{(N-4)/4} 
\eeq
for $N_1\ge N_2$, where $N=N_1+N_2$ and $\mu_1\ge \cdots\ge \mu_{N_2}$ are the eigenvalues of $\frac1{N_1}J^TJ$. 
The eigenvalues satisfy Condition~\ref{cond:empirical} and~\ref{cond:rigidity} with high probability; see Section~\ref{sec:proofeigval}.

\subsection{Critical point} \label{sec:cp}

We write 
\beq \label{eq:hatZ}
	\q= 
	- \int_{\gamma_1-\ii\infty}^{\gamma_1+\ii\infty} \int_{\gamma_2-\ii\infty}^{\gamma_2+\ii\infty} 
	e^{\n G(z_1, z_2)}  \, \dd z_2 \, \dd z_1
\eeq
with 
\beq
	G(z_1, z_2)= \ben (z_1+z_2) - \frac1{2\n} \sum_{i=1}^{\n} \log\left ( 4z_1z_2 - \mu_i(n) \right) - \aln \log z_1. 
\eeq

To evaluate the integral in \eqref{eq:hatZ} using the method of steepest-descent, 
we find the critical points of $G(z_1, z_2)$.
We have 
\beq \label{critical partial der}
	\partial_1 G= \ben - \frac{2z_2}{\n} \sum_{i=1}^{\n} \frac{1}{4z_1z_2 - \mu_i(n)} -\frac{\aln}{z_1}, \qquad 
		\partial_2 G = \ben - \frac{2z_1}{\n} \sum_{i=1}^{\n} \frac{1}{4z_1z_2 - \mu_i(n)} .
\eeq
Hence the critical points satisfy the equations
\beq \label{critical low 00temp1}
	z_1 - z_2 = \frac{\aln}{\ben}, \qquad \frac{z_1}{\n} \sum_{i=1}^{\n} \frac{1}{4z_1z_2 - \mu_i(n)} = \frac{\ben}{2}.
\eeq
Taking the imaginary parts, we find that at the critical points, 
\beqq
	\im z_1=\im z_2, 
	\qquad \frac1{\n} \sum_{i=1}^{\n} \frac{4|z_1|^2\im z_2 + \mu_i(n) \im z_1}{|4z_1z_2 - \mu_i(n)|^2}=0. 
\eeqq
Since $\mu_i(n)\ge 0$, $\im z_1=\im z_2=0$ at the critical points. 
Hence, all critical points, if exist, are real-valued.  

We now look for real critical points. 
Due to the branch cut of $G$, we look only for real critical points $(\gamma_1, \gamma_2)$ satisfying $4\gamma_1\gamma_2>\mu_1(n)$, $\gamma_1>0$, and $\gamma_2>0$. 
We set $4\gamma_1\gamma_2=\gamma$ and express the equations in terms of $\gamma_1$ and $\gamma$ instead of $\gamma_1$ and $\gamma_2$: 
\beq \label{critical low temp1}
	\gamma_1 - \frac{\gamma}{4\gamma_1} = \frac{\aln}{\ben}, \qquad 
	\frac{1}{\n} \sum_{i=1}^{\n} \frac{1}{\gamma - \mu_i(n)} = \frac{\ben}{2\gamma_1}
\eeq
where $\gamma>\mu_1(n)$. 
The first equation is a quadratic equation of $\gamma_1$ for given $\gamma$, and hence there are two solutions. 
Only one of them is positive given by 
\beq \label{eq:gamma_1 def}
\begin{split}
	\gamma_1 = \frac{\aln + \sqrt{\aln^2+\gamma \ben^2}}{2\ben} .
\end{split}
\eeq
This implies that 
\beq \label{eq:gamma_2 def}
\begin{split}
	\gamma_2 = \frac{-\aln + \sqrt{\aln^2+\gamma \ben^2}}{2\ben} .
\end{split}
\eeq
Inserting~\eqref{eq:gamma_1 def} into the second equation of~\eqref{critical low temp1}, we obtain an equation for $\gamma$ given by~\eqref{eq:gamma_def 1} below. 
The next lemma proves the existence and the uniqueness of the solution. 

\begin{lem} \label{lem:gamma_def}
The equation
\beq \label{eq:gamma_def 1}
	\frac{1}{\n} \sum_{i=1}^{\n} \frac{1}{\gamma - \mu_i(n)} 
	= \frac{\ben^2}{\aln + \sqrt{\aln^2+\gamma \be^2}}
\eeq
has a unique solution in the interval $(\mu_1(n), \infty)$.
\end{lem}

\begin{proof}
Let $L(\gamma)$ and $R(\gamma)$ be the left-hand side and right-hand side of \eqref{eq:gamma_def 1}, respectively. We observe that the function
$$
	f(\gamma)=\frac{c_1+\sqrt{c_2+\gamma}}{\gamma-\mu}
$$
has the derivative 
$$
	f'(\gamma)= \frac{-2c_1\sqrt{c_2+\gamma}- 2c_2 -\gamma-\mu}{2\sqrt{c_2+\gamma}(\gamma-\mu)^2}.
$$
Hence, if $c_1, c_2 > 0$ and $\mu>0$, then $h(\gamma)$ is a decreasing function of $\gamma \in (\mu, \infty)$. 
This shows that $\frac{L(\gamma)}{R(\gamma)}$ is a decreasing function of $\gamma\in (\mu_1(n), \infty)$. 
Since the equation \eqref{eq:gamma_def 1} is equivalent to $\frac{L(\gamma)}{R(\gamma)} = 1$, if the solution exists in the interval $(\mu_1(n), \infty)$, it is unique in the same interval.

We now prove the existence. We first notice that $L(\gamma) \to \infty$ as $\gamma \downarrow \mu_1(n)$ and $R(\gamma)$ is bounded above. Furthermore, $L(\gamma) = O(\gamma^{-1})$ as $\gamma \to \infty$ and $R(\gamma) \geq C \gamma^{-\frac{1}{2}}$ for some $C>0$ independent of $\gamma$. Thus,
\beq
	\lim_{\gamma \downarrow \mu_1} \frac{L(\gamma)}{R(\gamma)} = +\infty, \qquad \lim_{\gamma \to \infty} \frac{L(\gamma)}{R(\gamma)} = 0.
\eeq
Therefore, $\frac{L(\gamma)}{R(\gamma)} = 1$ has a unique solution in the interval $(\mu_1(n), \infty)$.
\end{proof}

In conclusion, 
\begin{enumerate}[(i)]
\item there are no critical values of $G$ with $\im z_1\neq 0$ or $\im z_2\neq 0$,  
\item there is a unique critical value $(\gamma_1, \gamma_2)$ such that $\gamma_1$ and $\gamma_2$ are real and positive, and $4\gamma_1\gamma_2>\mu_1(n)$,
\item the critical value $(\gamma_1, \gamma_2)$ is given by the formulas~\eqref{eq:gamma_1 def} and~\eqref{eq:gamma_2 def} where $\gamma\in (\mu_1(n), \infty)$ satisfies the equation~\eqref{eq:gamma_def 1}.
\end{enumerate}
Note that $(\gamma_1, \gamma_2)=(\gamma_1(n), \gamma_2(n))$ depends on $\n$ since $G$ depends on $\n$.

\subsection{Critical temperature} \label{subsec:critical temperature}

We discuss how we find the critical temperature formally from the critical point. 

Recall Condition~\ref{cond:aanndb} and Condition~\ref{cond:empirical}. 
Recall that $\red$ denotes the rightmost point of the support of $\wh\mu$. 
If $\gamma$ in~\eqref{eq:gamma_def 1} is $O(1)$ distance to the right of $\red$, then we may approximate the equation~\eqref{eq:gamma_def 1} by the $\n$-independent equation 
\beq \label{eq:ctemprinf}
	\int_{\R} \frac{1}{z-x} \dd \wh\mu(x)
	= \frac{\be^2}{\al + \sqrt{\al^2+z \be^2}}.
\eeq
Call the left-hand side and right-hand side by $L_\infty(z)$ and $R_\infty(z)$, respectively.
Note that $L_\infty(z)$ is well-defined for all real-valued $z\ge \red$ (and also non-real $z$). In particular, the integral converges when $z=\red$ due to Condition~\ref{cond:empirical}.
By the same calculation of the proof of Lemma~\ref{lem:gamma_def}, 
$\frac{L_\infty(z)}{R_\infty(z)}$ is a decreasing function of $z\in (\red, \infty)$. 
As before, $\frac{L_\infty(z)}{R_\infty(z)}\to 0$ as $z\to \infty$. 
However, unlike the previous lemma, the limit
\beq
	\lim_{z\downarrow \red} \frac{L_\infty(z)}{R_\infty(z)} = \frac{L_\infty(\red)}{R_\infty(\red)}.
\eeq
is finite. 
Hence the solution $z$ to the equation~\eqref{eq:ctemprinf} exists in $(\red,\infty)$ only if $L_\infty(\red)>R_\infty(\red)$, i.e., if 
\beq
	\int_{\R} \frac{1}{\red-x} \dd \wh\mu(x)
	> \frac{\sqrt{\al^2+\red \be^2}-\al}{\red}.
\eeq
Note that the left integral is a finite positive number due to the square-root vanishing assumption in Condition~\ref{cond:empirical}. 
Considered as a function of $\be$, the right-hand side $f(\be)$ is an increasing function of $\be$, $f(0)=0$, and $f(\be)\to +\infty$ as $\be\to \infty$. 
Hence the above inequality holds for all $\be<\be_c$ where $\be_c$ is defined by the equation
\beq \label{eq:B_c def}
	\int_{\R} \frac{1}{\red-x} \dd \wh\mu(x)
	= \frac{\sqrt{\al^2+\red \be_c^2}-\al}{\red}.
\eeq
Thus, we define the following critical value of $\be$.

\begin{defn} \label{defn:cv}
Define 
\beq \label{eq:B_c}
	\be_c= \sqrt{\red (\s(\red))^2+2\al\s(\red)} \quad \text{where } \; \s(z):= \int_{\R} \frac{1}{z-x} \dd \wh\mu(x).
\eeq
\end{defn}

The above discussion implies the following:
\begin{enumerate}[(a)]
\item For $0<\be < \be_c$, there is a unique solution $z_c$ in $(\red, \infty)$ to the equation~\eqref{eq:ctemprinf}. 
\item For $\be> \be_c $, there are no solutions to the equation~\eqref{eq:ctemprinf} in $(\red, \infty)$, and 
\end{enumerate}
We will show in Section~\ref{sec:high} that for the case (a), $\gamma$ in Lemma~\ref{lem:gamma_def} is indeed close to $z_c$. 
On the other hand, we will see in Section~\ref{sec:low} that for the case (b), the assumption that the point $\gamma$ in Lemma~\ref{lem:gamma_def} is $O(1)$ away from $\red$ is not true. 
This means that~\eqref{eq:ctemprinf} is not a good approximation to the equation~\eqref{eq:gamma_def 1}.

\subsection{Truncation of the double integral}

The following lemma gives an estimate on the double integral \eqref{eq:hatZ} outside a small disk of radius $N^{-\frac{1}{2}+\epsilon}$ about the point $(\gamma_1, \gamma_2)$. 
This result is used in later sections. 
The lemma does not require that $(\gamma_1, \gamma_2)$ is the critical point.

\begin{lem} \label{lem:estimate outside}
Let $\gamma_1=\gamma_1(\n)$, $\gamma_2=\gamma_2(\n)$ be any positive real numbers such that $4\gamma_1(\n)\gamma_2(\n)>\mu_1(\n)$ for all $\n$.
Suppose that there is a constant $C'>0$ such that 
$4\gamma_1(\n)\gamma_2(\n)-\mu_{\n}(\n)\le C'$ for all $\n$. 
Then, for any $\epsilon > 0$ and any $\Omega \subset \{ (y_1, y_2) \in \R^2 : y_1^2 + y_2^2 \geq \n^{-1+2\epsilon} \}$,
\beq
	 \int_\Omega \exp \left[ \n \re \left( G(\gamma_1+\ii y_1, \gamma_2+\ii y_2) - G(\gamma_1, \gamma_2) \right) \right]	\, \dd y_2 \, \dd y_1  \leq C e^{-\n^{\epsilon}}
\eeq
with high probability.
\end{lem}

\begin{proof}
We write $\mu_i(n)=\mu_i$ in this proof for a notational convenience. 
For $y_1, y_2 \in \R$, from the definition of $G$,
\beqq \begin{split}
	&\re \left[ G(\gamma_1+\ii y_1, \gamma_2+\ii y_2) - G(\gamma_1, \gamma_2) \right] \\
	&= - \frac{1}{4\n} \sum_{i=1}^{\n} \log \left[\left( 1- \frac{4y_1y_2}{4\gamma_1\gamma_2 - \mu_i} \right)^2 + 16 \left( \frac{\gamma_2 y_1 + \gamma_1 y_2}{4\gamma_1\gamma_2 - \mu_i} \right)^2 \right] 
	- \frac{\aln}{2} \log \left( 1 + \frac{y_1^2}{\gamma_1^2} \right).
\end{split} \eeqq

Consider the case $y_1y_2 \geq 0$. Then
\beqq \begin{split}
	&\left( 1- \frac{4y_1y_2}{4\gamma_1\gamma_2 - \mu_i} \right)^2 + 16 \left( \frac{\gamma_2 y_1 + \gamma_1 y_2}{4\gamma_1\gamma_2 - \mu_i} \right)^2 
	\geq 1- \frac{8y_1y_2}{4\gamma_1\gamma_2 - \mu_i} + 16 \left( \frac{\gamma_2 y_1 + \gamma_1 y_2}{4\gamma_1\gamma_2 - \mu_i} \right)^2 \\
	&\qquad \geq 1 + 8 \left( \frac{\gamma_2 y_1 + \gamma_1 y_2}{4\gamma_1\gamma_2 - \mu_i} \right)^2 \geq 1 + c(y_1^2 + y_2^2),
\end{split} \eeqq
where we used the fact that 
\beqq
	8 (\gamma_2 y_1 + \gamma_1 y_2)^2 \geq 32 \gamma_1 \gamma_2 y_1 y_2 \geq 8 y_1 y_2 (4\gamma_1\gamma_2 - \mu_i)
\eeqq
for the second inequality and that $|4\gamma_1\gamma_2 - \mu_i| < C$ uniformly for all $i$ in the third inequality.

For the case $y_1y_2 < 0$, we consider the following sub-cases:
\begin{enumerate}[(i)]
\item If $\gamma_2 |y_1| > 2 \gamma_1 |y_2|$, then $2|\gamma_2 y_1 + \gamma_1 y_2|\ge |\gamma_2 y_1|$, and hence
$$
	\left( 1- \frac{4y_1y_2}{4\gamma_1\gamma_2 - \mu_i} \right)^2 + 16 \left( \frac{\gamma_2 y_1 + \gamma_1 y_2}{4\gamma_1\gamma_2 - \mu_i} \right)^2 \geq 1 + 4 \left( \frac{\gamma_2 y_1}{4\gamma_1\gamma_2 - \mu_i} \right)^2 \geq 1 + c(y_1^2 + y_2^2).
$$

\item If $\gamma_2 |y_1| < \frac{1}{2} \gamma_1 |y_2|$,
then $2|\gamma_2 y_1 + \gamma_1 y_2|\ge |\gamma_1 y_2|$, and hence
$$
	\left( 1- \frac{4y_1y_2}{4\gamma_1\gamma_2 - \mu_i} \right)^2 + 16 \left( \frac{\gamma_2 y_1 + \gamma_1 y_2}{4\gamma_1\gamma_2 - \mu_i} \right)^2 \geq 1 + 4 \left( \frac{\gamma_1 y_2}{4\gamma_1\gamma_2 - \mu_i} \right)^2 \geq 1 + c(y_1^2 + y_2^2).
$$

\item If $\frac{1}{2} \gamma_1 |y_2| \leq \gamma_2 |y_1| \leq 2 \gamma_1 |y_2|$, then
$$
	\left( 1- \frac{4y_1y_2}{4\gamma_1\gamma_2 - \mu_i} \right)^2 + 16 \left( \frac{\gamma_2 y_1 + \gamma_1 y_2}{4\gamma_1\gamma_2 - \mu_i} \right)^2 \geq 1- \frac{8y_1y_2}{4\gamma_1\gamma_2 - \mu_i} \geq 1 + c(y_1^2 + y_2^2),
$$
since $-y_1y_2 =|y_1y_2| \geq c' (y_1^2 + y_2^2)$ for some $c'>0$.
\end{enumerate}

Thus, for all $y_1, y_2\in \R$, 
\beq \begin{split}
	&\re \left[ G(\gamma_1+\ii y_1, \gamma_2+\ii y_2) - G(\gamma_1, \gamma_2) \right] 
	\le - \frac{1}{4} \log \left( 1 + c(y_1^2 + y_2^2)\right). 
\end{split} \eeq
Now note that
\beq
	\log \left( 1 + c(y_1^2 + y_2^2)\right) 
	\ge \log \left( 1 + c\n^{-1+2\epsilon}\right) 
	\ge 	\frac{1}{2}c\n^{-1+2\epsilon} \quad \text{for $y_1^2+y_2^2\in [\n^{-1+2\epsilon}, \n]$}
\eeq
and
\beq
	\log \left( 1 + c(y_1^2 + y_2^2)\right) 
	\ge \log \left( c(y_1^2 + y_2^2) \right) \quad \text{for $y_1^2+y_2^2> \n$.} 
\eeq
Hence, 
\beq \begin{split}
	& \int_{\Omega} \exp \left[ \n \re \left( G(\gamma_1+\ii y_1, \gamma_2+\ii y_2) - G(\gamma_1, \gamma_2) \right) \right]	\, \dd y_2 \, \dd y_1  \\
	&\leq C e^{-\frac{c}{8}\n^{2\epsilon}} \int_{\n^{-1/2+\epsilon}}^{\n^{1/2}}  r \dd r + C \int_{\n^{1/2}}^{\infty}  (cr^2)^{-\n/4} r \dd r = O(e^{-c'\n^{2\epsilon}}) + O(\n^{-\n/4}).
\end{split} \eeq
This proves the lemma.
\end{proof}

\section{High temperature} \label{sec:high}

We consider the asymptotics of the double integral $\q$ in~\eqref{eq:hatZ} when $\be<\be_c$, 
where $\be_c$ is defined in~\eqref{eq:B_c}. We assume Conditions~\ref{cond:aanndb},~\ref{cond:empirical} and~\ref{cond:rigidity} throughout this section.
Recall that
\beq
	G(z_1, z_2)= \ben (z_1+z_2) - \frac1{2\n} \sum_{i=1}^{\n} \log\left ( 4z_1z_2 - \mu_i(n)\right) - \aln \log z_1. 
\eeq
As in our previous works \cite{BaikLee, BaikLee2}, we show that when $\be <\be_c$, the critical point of the random function $G$ is close to the critical point of a deterministic function.

Define
\beq \begin{split}
	G_\infty(z_1, z_2) = B (z_1+z_2) - \frac1{2} \int_{\R} \log\left ( 4z_1z_2 - x\right) \dd\wh\mu(x) - \al \log z_1 \,.
\end{split} \eeq
Then 
\beq \label{eq:Ginftygrad}
	\partial_1 G_\infty
	= B -  \int_{\R} \frac{2z_2}{4z_1z_2 - x} \dd\wh\mu(x)  - \frac{\al}{z_1} , 
	\qquad
	\partial_2 G_\infty
	= B - \int_{\R} \frac{2z_1}{4z_1z_2 - x} \dd\wh\mu(x).
\eeq
When $B<B_c$, the critical point of $G_\infty$ is given by
\beq \label{eq:z12ccp}
\begin{split}
	z_1^c = \frac{\al + \sqrt{\al^2+z_c \be^2}}{2\be},
	\qquad 
	z_2^c = \frac{-\al + \sqrt{\al^2+z_c \be^2}}{2\be},
\end{split}
\eeq
where $z_c$ is the solution to the equation
\beq \label{eq:ctemprinf2}
	\int_{\R} \frac{1}{z-x} \dd \wh\mu(x)
	= \frac{\be^2}{\al + \sqrt{\al^2+z \be^2}}
\eeq
satisfying $z_c\in (\red, \infty)$. 
We discussed in Subsection~\ref{subsec:critical temperature} that when $\be<\be_c$, there is unique $z_c$.

\bigskip

We start with the following lemma on the differences between the derivatives of $G$ and $G_\infty$, which is analogous to Lemma 5.1 of \cite{BaikLee}.

\begin{lem} \label{lem:G derivative high}
Fix $\theta > 0$ and set $B_{\theta} = \{ (z_1, z_2) : -\theta < \re (4z_1z_2) < \red +\theta, \, -\theta < \im (z_1z_2) < \theta \}$. Then the following hold. 
\begin{enumerate}[(i)]
\item For every $\epsilon > 0$ and each multi-index $m=(m_1, m_2)$, 
\beq
	\partial^m G(z_1, z_2) - \partial^m  G_{\infty}(z_1, z_2) = O(\n^{-1+\epsilon})
\eeq
uniformly on any compact subset of the region $\C^2 \setminus B_{\theta}$.

\item For every multi-index $m$,  $\partial^{m}G(z_1, z_2) = O(1)$ uniformly 
on any compact subset of $\C^2 \setminus B_{\theta}$.

\end{enumerate}
\end{lem}

\begin{proof}
(i) Let
\beq
	\wt G(z_1, z_2) = \be(z_1+z_2) - \frac1{2\n} \sum_{i=1}^{\n} \log\left ( 4z_1z_2 - g_i\right) - \al \log z_1
\eeq
where $g_i$ is the classical location of the $i$-th eigenvalue defined in \eqref{eq:classicallocationdef}. Then, from  Condition~\ref{cond:aanndb}, 
\beqq \begin{split}
	|G(z_1, z_2) - \wt G(z_1, z_2)| &= \left| (\ben-\be)(z_1+z_2)- \frac1{2\n} \sum_{i=1}^{\n} \log\left ( \frac{4z_1z_2 - \mu_i(n)}{4z_1z_2 - g_i} \right) -(\aln-\al)\log z\right| \\
	& \leq \frac1{2\n} \sum_{i=1}^{\n} \log \left( 1 + C|\mu_i - g_i| \right) + C' n^{-1-\delta} 
\end{split} \eeqq
uniformly on a compact subset of $\C^2 \setminus B_{\theta}$ since $|4z_1z_2-g_i|\ge c$. 
Hence, from the rigidity, Condition \ref{cond:rigidity},
\beqq
	|G(z_1, z_2) - \wt G(z_1, z_2)|  \leq \frac{C}{2\n} \sum_{i=1}^{\n} |\mu_i - g_i|  + \frac{C'}{n^{1+\delta}}
	\leq \frac{C\n^{\epsilon}}{\n}
\eeqq
in any compact subset of $\C^2 \setminus B_{\theta}$.
We now compare $\wt G(z_1, z_2)$ with $G_{\infty}(z_1, z_2)$. 
For $2\le i\le \n-1$,
\beqq
	\int_{g_i}^{g_{i-1}} \log\left ( 4z_1z_2 - x\right) \dd\wh\mu(x)
	 \leq \frac{1}{\n} \log\left ( 4z_1z_2 - g_i \right) \leq \int_{g_{i+1}}^{g_i} \log\left ( 4z_1z_2 - x\right) \dd\wh\mu(x).
\eeqq
Summing over the index $i$ and using the trivial estimates 
$$
	\int_{g_2}^{\red} \log(4z_1z_2 -x) \dd\wh\mu(x) = O(\n^{-1}), \qquad \int_{0}^{g_{\n-2}} \log(4z_1z_2 -x) \dd\wh\mu(x) = O(\n^{-1}),
$$
and $\log (4z_1z_2 -g_i) = O(1)$ for any compact subset of $\C^2 \setminus B_{\theta}$, we find that $\wt G(z_1, z_2) - G_{\infty}(z_1, z_2) = O(\n^{-1})$. 
Hence, $G(z_1, z_2) - G_{\infty}(z_1, z_2) = O(\n^{-1+\epsilon})$. 
For the derivatives, the function $ \log\left ( 4z_1z_2 - x\right)$ is replaced by $\frac1{ ( 4z_1z_2 - x)^k}$ for positive integers $k$, and the proof is almost same. 
  
(ii) can be proved in a similar manner since, for any compact subset of $\C^2 \setminus B_{\theta}$, $\log (4z_1z_2 -g_i) = O(1)$ and $\frac1{ ( 4z_1z_2 - x)^k}=O(1)$ for positive integers $k$.
\end{proof}

We now compare the critical point $(\gamma_1, \gamma_2)$ of $G$ and the critical point $(z_1^c, z_2^c)$ of $G_\infty$. 
Recall that $(\gamma_1, \gamma_2)$ depends on $\n$. 

\begin{lem} \label{lem:G critical estimate}
For any $\epsilon > 0$, the following hold.
\begin{enumerate}[(i)]
\item We have 
\beq \label{eq:gamma12 estimate}
	\gamma_1 - z_1^c = O(\n^{-1+\epsilon}), \qquad  \gamma_2 - z_2^c = O(\n^{-1+\epsilon}).
\eeq
\item There is a  positive constant $c$, independent of $\n$, such that 
\beq
	\text{$4\gamma_1 \gamma_2 - \mu_1(n) >c$ and $4\gamma_1 \gamma_2 - \red >c$.}
\eeq 
\item We have 
\beq
	G(\gamma_1, \gamma_2) = G(z_1^c, z_2^c) + O(\n^{-2+\epsilon})
\eeq
and for any multi-index $m=(m_1, m_2)$ satisfying $|m|>0$, 
\beq
	\partial^{m} G(\gamma_1, \gamma_2) = \partial^{m} G(z_1^c, z_2^c) + O(\n^{-1+\epsilon}).
\eeq
\end{enumerate}
\end{lem}

\begin{proof}
(i) We first show that $\gamma=4\gamma_1\gamma_2$ and $z_c=4z_1^cz_2^c$ satisfy 
$\gamma - z_c = O(n^{-1+\epsilon})$.
The value $\gamma$ is determined by the equation in~\eqref{eq:gamma_def 1}, which can be written as 
$L(x)=R(x)$ where 
\beqq
	L(x):= \frac{1}{\n} \sum_{i=1}^{\n} \frac{1}{x - \mu_i(n)} , 
	\qquad 
	R(x):= \frac{\ben^2}{\aln + \sqrt{\aln^2+x \ben^2}} .
\eeqq
Similarly, the point $z_c$ is a solution of the equation $L_\infty(x)=R(x)$ where  
\beqq 
	L_\infty(x):= \int_{\R} \frac{\dd \wh \mu(y)}{x- y}, 
	\qquad R_\infty(x):= \frac{\be^2}{\al + \sqrt{\al^2+x \be^2}}.
\eeqq
We showed in the proof of Lemma \ref{lem:gamma_def} that $F(x) := \frac{R(x)}{L(x)}$ satisfies
$F'(x)<0$ for all $x>\mu_1(n)$. 
The same calculation shows that $F_\infty(x):=\frac{R_\infty(x)}{L_\infty(x)}$ satisfies $F_\infty'(x)<0$ for all $x>\red$.
Since
\beqq
	L(x)= \frac12(\ben-\partial_2 G(1, x)), 
	\qquad
	L_\infty(x)= \frac12(\be-\partial_2 G_\infty(1, x)),	
\eeqq
we find from Lemma~\ref{lem:G derivative high} (i) that $F(x)=F_\infty(x)+O(\n^{-1+\frac{\epsilon}2})$ uniformly for $x$ in any compact subset of the interval $(\red, \infty)$. 
Note that we used $\epsilon/2$ when we apply Lemma~\ref{lem:G derivative high}.
Recall that $z_c>\red$. 
Hence, $F(z_c)=F_\infty(z_c)+O(N^{-1+\frac{\epsilon}2})= 1+ O(\n^{-1+\frac{\epsilon}2})$.
By Taylor series, 
\beqq \begin{split}
	F(z_c\pm \n^{-1+\epsilon})
	&= F_\infty(z_c\pm \n^{-1+\epsilon})+O(\n^{-1+\frac{\epsilon}2}) = 1 \pm F_\infty'(z_c) \n^{-1+\epsilon}+O(\n^{-2+2\epsilon})+O(\n^{-1+\frac{\epsilon}2}). 
\end{split} \eeqq
Since $F_\infty'(z_c)<0$, we find that 
\beqq
	F(z_c- \n^{-1+\epsilon})<1, \qquad F(z_c- \n^{-1+\epsilon})>1. 
\eeqq
This implies that 
\beqq
	\gamma\in (z_c-\n^{-1+\epsilon}, z_c+\n^{-1+\epsilon}).
\eeqq
We obtain (i) since $\gamma_1$ and $\gamma_2$ are given in terms of $\gamma$ by~\eqref{eq:gamma_1 def} and~\eqref{eq:gamma_2 def}, and respectively, and $z_1^c$ and $z_2^c$ are given by~\eqref{eq:z12ccp} in terms of $z_c$. 

(ii) follows from (i).

(iii) From the Taylor expansion and the bounds in Lemma \ref{lem:G derivative high} (ii), 
\beqq \begin{split}
	G(z_1^c, z_2^c) &= G(\gamma_1, \gamma_2) + \partial_1 G(\gamma_1, \gamma_2) (z_1^c-\gamma_1) + \partial_2 G(\gamma_1, \gamma_2) (z_2^c-\gamma_2) + O(\n^{-2+\epsilon}) \\
	&= G(\gamma_1, \gamma_2) + O(\n^{-2+\epsilon})
\end{split} \eeqq
and
\beqq
	\partial^{m} G(z_1^c, z_2^c) = \partial^{m} G(\gamma_1, \gamma_2) + O(\n^{-1+\epsilon})
\eeqq
for any multi-index $m$ satisfying $|m|>0$. This completes the proof of the lemma.
\end{proof}

We evaluate the integral~\eqref{eq:hatZ} using the method of steepest-descent.

\begin{lem} \label{lem:integral high}
Let $\be < \be_c$ for $\be_c$ defined in~\eqref{eq:B_c}. Then for any $\epsilon > 0$,
\beq
	\q
	= e^{\n G(\gamma_1, \gamma_2)} \frac{\pi}{\n \sqrt{D(\gamma_1, \gamma_2)}} \left( 1 + O(\n^{-1+\epsilon}) \right)
\eeq
where $D(\gamma_1, \gamma_2)$ is the discriminant 
$$
	D(\gamma_1, \gamma_2) = \partial_1^2 G(\gamma_1, \gamma_2) \cdot \partial_2^2 G(\gamma_1, \gamma_2) - \big( \partial_1 \partial_2 G(\gamma_1, \gamma_2) \big)^2.
$$
\end{lem}

\begin{proof}
Changing the variables, 
\beqq \begin{split} 
	\q	&= \frac{1}{\n} e^{\n G(\gamma_1, \gamma_2)} \int_{-\infty}^{\infty} \int_{-\infty}^{\infty} \exp \left[ \n \left( G(\gamma_1+\ii \frac{t_1}{\sqrt{\n}}, \gamma_2+\ii \frac{t_2}{\sqrt {\n}}) - G(\gamma_1, \gamma_2) \right) \right]	\, \dd t_2 \, \dd t_1 .
\end{split} \eeqq
Lemma \ref{lem:estimate outside} shows that, the part of the last double integral over the region $\R^2 \setminus [-\n^{\epsilon}, \n^{\epsilon}]^2$ is $O(e^{-\n^\epsilon})$.
On the other hand, for $|t_1|, |t_2| \leq {\n}^{\epsilon}$, 
\beqq \begin{split}
	&G(\gamma_1+\ii \frac{t_1}{\sqrt{\n}}, \gamma_2+\ii \frac{t_2}{\sqrt{\n}}) - G(\gamma_1, \gamma_2) \\
	&= -\frac{1}{2\n} \left( \partial_1^2 G(\gamma_1, \gamma_2) t_1^2 + 2\partial_1 \partial_2 G(\gamma_1, \gamma_2) t_1 t_2 + \partial_2^2 G(\gamma_1, \gamma_2) t_2^2 \right) \\
	& \quad -\frac{\ii}{6\n^{\frac{3}{2}}} \left( \partial_1^3 G(\gamma_1, \gamma_2) t_1^3 +3\partial_1^2 \partial_2 G(\gamma_1, \gamma_2) t_1^2 t_2 + 3\partial_1 \partial_2^2 G(\gamma_1, \gamma_2) t_1 t_2^2 + \partial_2^3 G(\gamma_1, \gamma_2) t_2^3 \right) +O(\n^{-2+4\epsilon}) \\
	& =: -\frac{X_2 (t_1, t_2)}{\n} - \frac{\ii X_3 (t_1, t_2)}{\n^{\frac{3}{2}}} + O(\n^{-2+4\epsilon})
\end{split} \eeqq
where we used Lemma~\ref{lem:G critical estimate} (i) and Lemma~\ref{lem:G derivative high} (ii) for the error estimate. 
Hence,
\beqq \begin{split}
	&\int_{-\n^{\epsilon}}^{\n^{\epsilon}} \int_{-\n^{\epsilon}}^{\n^{\epsilon}} \exp \left[ \n \left( G(\gamma_1+\ii \frac{t_1}{\sqrt{\n}}, \gamma_2+\ii \frac{t_2}{\sqrt{\n}}) - G(\gamma_1, \gamma_2) \right) \right]	\, \dd t_2 \, \dd t_1 \\
	&=\int_{-\n^{\epsilon}}^{\n^{\epsilon}} \int_{-\n^{\epsilon}}^{\n^{\epsilon}} e^{-X_2 (t_1, t_2)} \dd t_2 \, \dd t_1 - \ii \int_{-\n^{\epsilon}}^{\n^{\epsilon}} \int_{-\n^{\epsilon}}^{\n^{\epsilon}} \frac{X_3 (t_1, t_2)}{\sqrt{\n}} e^{-X_2 (t_1, t_2)} \dd t_2 \, \dd t_1 + O(\n^{-1+6\epsilon}).
\end{split} \eeqq
Since 
\beqq
X_3 (-t_1, -t_2) e^{-X_2 (-t_1, -t_2)} = -X_3 (t_1, t_2) e^{-X_2 (t_1, t_2)},  
\eeqq
the integral in the middle vanishes. 
On the other hand, from the estimate $\int_{\n^{\epsilon}}^{\infty} e^{-t^2} \dd t = O(\n^{-\epsilon} e^{-\n^{2\epsilon}})$, we obtain that
\beqq \begin{split}
	\int_{-\n^{\epsilon}}^{\n^{\epsilon}} \int_{-\n^{\epsilon}}^{\n^{\epsilon}} e^{-X_2 (t_1, t_2)} \dd t_2 \, \dd t_1
	= \frac{\pi}{\sqrt{D(\gamma_1, \gamma_2)}} + O(\n^{-1+6\epsilon}).
\end{split} \eeqq
Thus, we obtain the lemma. 
\end{proof}

The following is the main result for the double integral $\q$ when $\be<\be_c$.

\begin{prop}[Random double integral for high temperature] \label{prop:integral high}
Assume Conditions~\ref{cond:aanndb},~\ref{cond:empirical} and~\ref{cond:rigidity}. 
Define
\beq
	\h(z):= \int_{\R} \log(z-x) \dd \wh\mu(x)\,.
\eeq
Suppose that $\be$ in Condition~\ref{cond:aanndb} satisfies $0<\be<\be_c$ where $\be_c$ is defined in \eqref{eq:B_c}. 
Then, setting $z_c$ be the unique solution of the equation~\eqref{eq:ctemprinf}, 
\beq \label{eq:cpeq}
	\h'(z_c) = \frac{\be^2}{\al + \sqrt{\al^2+z_c \be^2}}, \qquad z_c\in (\red, \infty), 
\eeq
we have for every $\epsilon > 0$,
\beq
	\frac1{\n}\log \q
	= \aw - \frac1{2\n} \left[  \sum_{i=1}^{\n} \log\left ( z_c - \mu_i\right) - \n \h(z_c) \right] - \frac{\log \n}{\n} + \frac1{2\n} \log \left( \frac{\pi^2}{\dw}  \right) + O(\n^{-2+\epsilon}) 
\eeq
where 
\beq \begin{split}
	\aw &= \sqrt{\al^2+z_c \be^2} - \al \log \left( \frac{\al + \sqrt{\al^2+z_c \be^2}}{2\be} \right)  - \frac1{2} \h(z_c),  \\
	\dw &=  -8 \al \h''(z_c) -8z_c \h'(z_c)\h''(z_c) -4 (\h'(z_c))^2 .
\end{split} \eeq
\end{prop}

\begin{proof}
From Lemma~\ref{lem:integral high}, 
\beqq \begin{split}
	\frac1{\n}\log \q
	&=G(\gamma_1, \gamma_2) - \frac{\log \n}{\n} + \frac1{2\n} \log \left( \frac{\pi^2}{D(\gamma_1, \gamma_2)} \right) + O(\n^{-2+\epsilon}) .
\end{split} \eeqq
Using Lemma~\ref{lem:G critical estimate} (iii), we write
\beqq
	G(\gamma_1, \gamma_2) = G(z_1^c, z_2^c) + O(\n^{-2+\epsilon})
	= G_\infty(z_1^c, z_2^c) + \left[ G(z_1^c, z_2^c) - G_\infty(z_1^c, z_2^c) \right] + O(\n^{-2+\epsilon}) .
\eeqq
We have
\beqq
	G(z_1^c, z_2^c)- G_\infty(z_1^c, z_2^c) = -\frac1{2\n} \left[  \sum_{i=1}^{\n} \log\left ( z_c - \mu_i\right)
	-  \n \int_{\R} \log\left ( z_c - x\right) \dd\wh\mu(x) \right] + O(\n^{-1-\delta}) .
\eeqq
We also have 
\beqq \begin{split}
	G_\infty(z_1^c, z_2^c) = \sqrt{\al^2+z_c \be^2} - \frac1{2} \int_{\R} \log\left ( z_c - x\right) \dd\wh\mu(x) - \al \log \left( \frac{\al + \sqrt{\al^2+z_c \be^2}}{2\be} \right)= \aw \,.
\end{split} \eeqq
It remains to compare $D(\gamma_1, \gamma_2)$ with $\dw$. Using Lemma~\ref{lem:G critical estimate} (iii) and Lemma~\ref{lem:G derivative high} (i), 
\beqq
	D(\gamma_1, \gamma_2) =D_\infty (z_1^c, z_2^c) +  O(\n^{-1+\epsilon})
\eeqq
where
\beqq
	D_\infty (z_1^c, z_2^c) 
	:= \partial_1^2 G_\infty(z_1^c, z_2^c)  \cdot \partial_2^2 G_\infty(z_1^c, z_2^c) - \big( \partial_1 \partial_2 G_\infty(z_1^c, z_2^c) \big)^2 .
\eeqq
From direct computation, 
\beqq
	D_\infty (z_1^c, z_2^c) 
	=  -8 \al \h''(z_c) -8z_c \h'(z_c)\h''(z_c) -4 (\h'(z_c))^2= \dw.
\eeqq
This completes the proof.
\end{proof}

\section{Low temperature} \label{sec:low}

In this section, we consider the asymptotics of 
\beq \begin{split} \label{double integral estimate 1}
	e^{-\n G(\gamma_1, \gamma_2)} \q	= - \int_{\gamma_1-\ii\infty}^{\gamma_1+\ii\infty} \int_{\gamma_2-\ii\infty}^{\gamma_2+\ii\infty} \exp \left[ \n \left( G(z_1, z_2) - G(\gamma_1, \gamma_2) \right) \right] \, \dd z_2 \, \dd z_1 \,.
\end{split} \eeq
when $\be> \be_c$. We assume Conditions~\ref{cond:aanndb},~\ref{cond:empirical} and~\ref{cond:rigidity} throughout this section.

Unlike the previous section, when $\be>\be_c$, the critical point $(\gamma_1, \gamma_2)$ of $G$ is not approximated by the critical point of $G_{\infty}$. Indeed, we showed in Subsections~\ref{sec:cp} and~\ref{subsec:critical temperature} that $G_\infty$ has no critical point when $\be>\be_c$, while $(\gamma_1, \gamma_2)$ exists for all $\be$.
We show in Lemma~\ref{lem:lowcp} below that $\gamma=4\gamma_1\gamma_2$ is actually close to the branch point $\mu_1(\n)$.  
Due to this fact, the control of the double integral becomes subtle.  
We had a similar situation for a random single integral in \cite{BaikLee} for the usual SSK model. 
In this paper, we have a double integral, and this brings an additional difficulty. In particular, the symmetry we used in \cite{BaikLee}, which simplified the analysis, is no longer valid. 
In the below, we will choose the integration contours in a certain explicit way and show that it is possible to reduce the double integral to the product of two single integrals plus an error. 
One of the single integral is trivial and the other single integral has a certain symmetry that can be used to simplify the method  of steepest-descent in a manner similar to the analysis of \cite{BaikLee}.

In Subsections~\ref{subsec:gamma a priori}--\ref{subsec:I_12}, we prove the following lemma.
The conclusion of this section is given in Subsection~\ref{sec:lowtempresult}. 

\begin{lem} \label{lem:lowt}
Assume Conditions~\ref{cond:aanndb},~\ref{cond:empirical} and~\ref{cond:rigidity}.
Suppose that $\be$ in Condition~\ref{cond:aanndb} satisfies $\be > \be_c$.
Let $(\gamma_1, \gamma_2)$ be the critical point of $G$ given by~\eqref{eq:gamma_1 def} and~\eqref{eq:gamma_2 def}. Then, for every $\epsilon>0$, there is a constant $C>0$ such that 
\beq
	C\n^{-\frac32-\epsilon}\le e^{-\n G(\gamma_1, \gamma_2)}  \q \le C\n^{-\frac12}.
\eeq
\end{lem}

\begin{rem}[Notational Remark]
In order to lighten up the notations, we will  write $\mu_i$ for $\mu_i(\n)$ in the rest of this section. 
It should be understood that $\mu_i$ depends on $\n$. 
\end{rem}

\subsection{A priori estimate on $\gamma$} \label{subsec:gamma a priori}

We begin by approximating $\gamma = 4\gamma_1\gamma_2$ and introducing a priori estimates that will be used in this section.

\begin{lem} \label{lem:lowcp}
For any $0 < \epsilon < 1$, the solution $\gamma$ in Lemma \ref{lem:gamma_def},
\beq \label{eq:gmen}
	\frac{1}{\n} \sum_{i=1}^{\n} \frac{1}{\gamma - \mu_i} 
	= \frac{\ben^2}{\aln + \sqrt{\aln^2+\gamma \ben^2}}, \qquad \gamma>\mu_1,
\eeq
satisfies the inequality 
\beq \label{eq:gammalocationesm}
	\frac{\red^{1/2}}{2\be\n} \leq \gamma - \mu_1 \leq \frac{n^{\epsilon}}{n}.
\eeq	
\end{lem}

\begin{proof}
We follow the proof of Lemma 6.1 in \cite{BaikLee}. Define $L(\gamma)$ and $R(\gamma)$ to be the left-hand side and the right-hand side of \eqref{eq:gmen}, respectively. The equation~\eqref{eq:gmen} is equivalent to the equation $\frac{L(\gamma)}{R(\gamma)}=1$. 
Since $\mu_1\to \red$ and $\ben\to \be$,  
\beqq
	\mu_1 + \frac{\red^{1/2}}{2\be\n} \geq \frac{\red}{2}. 
\eeqq
Since $L(\gamma)\ge \frac1{n(\gamma-\mu_1)}$ for $\gamma>\mu_1$, we find that 
\beq
	L\big(\mu_1 + \frac{\red^{1/2}}{2\be\n} \big) \geq 
	\frac{2\be}{\red^{1/2}} \geq \frac{\ben^2}{\sqrt{ \big(\mu_1 + \frac{\red^{1/2}}{2\be\n}\big) \ben^2}} \geq R\big(\mu_1 + \frac{\red^{1/2}}{2\be\n} \big).
\eeq
Since $L(x)/R(x)$ is a decreasing function of $x$ (see the proof of Lemma~\ref{lem:gamma_def}), 
this implies the lower bound of~\eqref{eq:gammalocationesm}.

The upper bound is proved if we show that $L(\mu_1 + \n^{-1+4\epsilon}) < R(\mu_1 + \n^{-1+4\epsilon})$ for any $0< \epsilon < \frac{1}{4}$.
From Condition \ref{cond:rigidity}, $|\mu_i - g_i| \leq \n^{-2/3}$ for $\n^{3\epsilon} \leq i \leq \n -\n^{3\epsilon}$. For such $i$, we note that $\red - g_i \geq c\n^{-2/3 +2\epsilon}$. Since $\mu_1 = \red +O(\n^{-2/3 +\epsilon})$,
\beqq
	\frac{1}{n} \sum_{i=\n^{3\epsilon}}^{\n-\n^{3\epsilon}} \frac{1}{\mu_1 + \n^{-1+4\epsilon} -\mu_i} = \frac{1}{n} \sum_{i=\n^{3\epsilon}}^{\n-\n^{3\epsilon}} \frac{1}{\red -g_i} \left( 1+ O(\n^{-\epsilon}) \right).
\eeqq
Approximating the last sum by an integral as in the proof of Lemma \ref{lem:G derivative high}, we find that
\beqq
	\left| \frac{1}{n} \sum_{i=\n^{3\epsilon}}^{\n-\n^{3\epsilon}} \frac{1}{\mu_1 + \n^{-1+4\epsilon} -\mu_i} - \int_{\led}^{\red} \frac{\dd\wh\mu (x)}{\red -x} \right| = O(\n^{-1/3 +\epsilon}).
\eeqq
(See also Equations (6.6) and (6.7) in \cite{BaikLee}.)
For $1 \leq i < \n^{3\epsilon}$, since $\mu_1 \geq \mu_i$,
\beqq
	\frac{1}{n} \sum_{i=1}^{\n^{3\epsilon}-1} \frac{1}{\mu_1 + \n^{-1+4\epsilon} -\mu_i} = O(\n^{-\epsilon}).
\eeqq
Finally, for $\n-\n^{3\epsilon} < i \leq \n$, since $\mu_1 - \mu_i > c > 0$,
\beqq
	\frac{1}{n} \sum_{i=\n-\n^{3\epsilon}+1}^{\n} \frac{1}{\mu_1 + \n^{-1+4\epsilon} -\mu_i} = O(\n^{-1+3\epsilon}).
\eeqq
Combining the estimates, we find that
\beq
	L(\mu_1 + \n^{-1+4\epsilon}) = \int_{\led}^{\red} \frac{\dd\wh\mu (x)}{\red -x} + O(\n^{-\epsilon}).
\eeq
On the other hand, since $\mu_1 + \n^{-1+4\epsilon} \to \red$, 
\beq
	R(\mu_1 + \n^{-1+4\epsilon}) 
	\to  \frac{\be^2}{\al + \sqrt{\al^2+\red \be^2}} 
	= \frac{\sqrt{\al^2+\red \be^2}-\al}{\red} .
\eeq
From the definition of $\be_c$ in \eqref{eq:B_c def},
\beqq
	\frac{\sqrt{\al^2+\red \be^2}-\al}{\red}
	> \frac{\sqrt{\al^2+\red \be_c^2}-\al}{\red} = \int_{\led}^{\red} \frac{\dd\wh\mu (x)}{\red -x} .
\eeqq
Hence 
\beq
	R(\mu_1 + \n^{-1+4\epsilon})  > L(\mu_1 + \n^{-1+4\epsilon}) +c
\eeq
for some $c>0$ for all large enough $n$. 
This proves the lemma.
\end{proof}

Since $\gamma$ is well approximaed by $\mu_1$ and $\mu_1$ is close to $\red$, heuristically,
\beq
	\frac{1}{\n} \sum_{i=1}^{\n} \log(\gamma - \mu_i) \approx \frac{1}{\n} \sum_{i=1}^{\n} \left[ \log(\red - \mu_i) + \frac{\gamma - \red}{\red - \mu_i} \right] \approx \h(\red) + (\mu_1 - \red) \h'(\red).
\eeq
In the following lemma, we describe the approximation above rigorously and also estimate $\sum_{i=1}^{\n} (\gamma - \mu_i)^{-\ell}$ for $\ell = 2, 3, \dots$. Since the following lemma can be proved in a similar manner to the proof of Lemma 6.2 of \cite{BaikLee}, we omit the proof.

\begin{lem} \label{lem:logsummuigamma}
Recall the definition of $\h(z)$ in Proposition \ref{lem:integral high}. Then, for any $0 < \epsilon < 1$,
\beq
	\frac{1}{\n} \sum_{i=1}^{\n} \log(\gamma - \mu_i) =  \h(\red)+(\mu_1-\red)\h'(\red) + O(\n^{-1+\epsilon}).
\eeq
Furthermore, for any $0 < \epsilon < 1$ there is a constant $C_0 > 0$ such that
\beq \label{G^n bound}
	\n^{\ell (1-\epsilon)} \leq \sum_{i=1}^{\n} \frac{1}{(\gamma - \mu_i)^{\ell}} \leq C_0^{\ell} \n^{\ell+\epsilon}
\eeq
for all $\ell = 2, 3, \dots$. Here, $C_0$ does not depend on $\ell$.
\end{lem}

\begin{proof}
See Lemma 6.2 of \cite{BaikLee}.
\end{proof}

\subsection{Truncation and deformation of the coutour} \label{subsec:truncate}

In Subsections~\ref{subsec:truncate}--\ref{subsec:I_12}, we fix $0 < \epsilon < \frac{1}{100}$ and prove Lemma~\ref{lem:lowt}. 

Lemma \ref{lem:estimate outside} implies that the part of the double integral~\eqref{double integral estimate 1} with $|\im z_1| \geq \n^{-\frac{1}{2}+\epsilon}$ is $O(e^{-\n^\epsilon})$. 

For the part $|\im z_1| < \n^{-\frac{1}{2}+\epsilon}$, we deform the $z_2$-integral to a different vertical contour passing through a new point $\wt\gamma_2$ such that the difference $|G(\gamma_1, \gamma_2) - G(\gamma_1 + \ii y_1, \wt\gamma_2)|$ is sufficiently small. Intuitively, since the main contribution to the change of $G(z_1, z_2)$ near the critical point comes from the term $\frac{1}{4z_1z_2 -\mu_1}$, it must be very sensitive to the change of the product $z_1z_2$ but not to the change of the individual variable $z_1$ or $z_2$ while $z_1z_2$ is fixed. Thus, for $y_1\in \R$, we define
\beq
	\wt \gamma_2 \equiv \wt \gamma_2(y_1) = \frac{\gamma_1 \gamma_2}{\gamma_1 + \ii y_1}
\eeq
and analyze the double integral in \eqref{double integral estimate 1} with the deformed contour that passes through $\wt\gamma_2$ for the $z_2$-integral.

Before we peform the analysis, we check that it is possible to deform the contour $\gamma_2+\ii \R$ to $\wt \gamma_2+\ii \R$ for given $z_1 \in  \gamma_1 + \ii \R$.  
For fixed $z_1 = \gamma_1 + \ii y_1$, the branch cut $\Gamma_c$ of the logarithmic function in $G(z_1, z_2)$ as a function of $z_2$ is
$$
	\Gamma_c = \{ z_2 \in \C: 4z_1z_2-\mu_1 \in \R^{-} \cup \{ 0 \} \}.
$$
If $z_2 \in \Gamma_c$, then there exists $r \geq 0$ such that
\beqq
	z_2 = \frac{\mu_1 -r}{4(\gamma_1 +\ii y_1)} = \frac{\mu_1 -r}{4\gamma_1 \gamma_2} \wt \gamma_2. 
\eeqq
Since $4\gamma_1 \gamma_2 > \mu_1$, this implies that $\re z_2 < \re \wt \gamma_2$, and hence $\Gamma_c$ does not intersect the half plane $\{ z \in \C : \re z\ge \re \wt \gamma_2 \}$.
Therefore, we can deform the  $z_2$-contour, and hence 
\beq \begin{split} \label{double integral estimate 2}
	&-\int_{\gamma_1-\ii \n^{-\frac{1}{2}+\epsilon}}^{\gamma_1+\ii \n^{-\frac{1}{2}+\epsilon}} \int_{\gamma_2-\ii\infty}^{\gamma_2+\ii\infty} \exp \left[ \n \left( G(z_1, z_2) - G(\gamma_1, \gamma_2) \right) \right] \, \dd z_2 \, \dd z_1 \\
	&\quad =  \int_{-\n^{-\frac{1}{2}+\epsilon}}^{\n^{-\frac{1}{2}+\epsilon}} \int_{-\infty}^\infty  \exp \left[ \n \left( G(\gamma_1+\ii y_1, \wt \gamma_2+\ii y_2) - G(\gamma_1, \gamma_2) \right) \right] \, \dd y_2 \, \dd y_1 \,.
\end{split} \eeq
Recall that $\wt \gamma_2 \equiv \wt \gamma_2(y_1)$ depends on $y_1$.

We now truncate the $y_2$-integral. 
From the definition of $\wt \gamma_2$ and $G$, for all $y_1, y_2\in\R$, 
\beq \label{eq:G diff 0101} \begin{split}
	&G(\gamma_1 + \ii y_1, \wt \gamma_2 + \ii y_2) - G(\gamma_1, \gamma_2) 
	= \ii \ben (y_1 + y_2) -\ii \ben \frac{\gamma_2 y_1}{\gamma_1+\ii y_1}  \\
	&\qquad\quad  - \aln \log \left( 1+\frac{\ii y_1}{\gamma_1} \right) - \frac{1}{2\n} \sum_{i=1}^{\n} \log\left( 1 -\frac{4y_1 y_2}{4\gamma_1 \gamma_2 - \mu_i} + \frac{4\ii \gamma_1 y_2}{4\gamma_1 \gamma_2 - \mu_i} \right) .
\end{split} \eeq

\begin{lem} \label{lem:y_2 cutoff00}
Uniformly for $|y_1| \leq \n^{-\frac{1}{2}+\epsilon}$,
\beq \label{eq:y2cuff00}
	\left( \int_{-\infty}^{-\n^{-\frac12+2\epsilon}} + \int_{\n^{-\frac12+2\epsilon}}^\infty \right) \exp \left[ \n \left( G(\gamma_1+\ii y_1, \wt \gamma_2+\ii y_2) - G(\gamma_1, \gamma_2) \right) \right] \, \dd y_2
	= O(e^{-\n^{\epsilon}})
\eeq
\end{lem}

\begin{proof}
The proof is similar to Lemma \ref{lem:estimate outside}, but easier. 
Taking the real part of~\eqref{eq:G diff 0101}, 
\beq \label{eq:Gdiffincutf} \begin{split}
	& \re \left( G(\gamma_1+\ii y_1, \wt \gamma_2+\ii y_2) - G(\gamma_1, \gamma_2) \right) \\
	& = -\frac{\ben \gamma_2y_1^2}{\gamma_1^2+y_1^2}
	- \frac{\aln}{2} \log \left( 1+ \frac{y_1^2}{\gamma_1^2} \right) 
	- \frac1{4\n} \sum_{i=1}^{\n} \log \left[ \left( 1- \frac{4y_1y_2}{4\gamma_1 \gamma_2 -\mu_1}\right)^2 + \left( \frac{4\gamma_1 y_2}{4\gamma_1 \gamma_2 -\mu_i}  \right)^2 \right] . 
\end{split}
\eeq
If $y_1y_2\le 0$, then 
\beqq
	\left( 1- \frac{4y_1y_2}{4\gamma_1 \gamma_2 -\mu_1}\right)^2 + \left( \frac{4\gamma_1 y_2}{4\gamma_1 \gamma_2 -\mu_i}  \right)^2 
	\ge 1+  \left( \frac{4\gamma_1 y_2}{4\gamma_1 \gamma_2 -\mu_i}  \right)^2 \ge 1+cy_2^2 .
\eeqq
If $y_1y_2\ge 0$, then 
\beqq
	\left( 1- \frac{4y_1y_2}{4\gamma_1 \gamma_2 -\mu_1}\right)^2 + \left( \frac{4\gamma_1 y_2}{4\gamma_1 \gamma_2 -\mu_i}  \right)^2 
	\ge 1 - \frac{8 y_1y_2}{4\gamma_1 \gamma_2 -\mu_1}+  \frac{16\gamma_1^2y_2^2}{(4\gamma_1 \gamma_2 -\mu_i)^2} 
	\ge 1+cy_2^2 .
\eeqq
since $y_1y_2 \le y_2^2 \n^{-\epsilon}$ for $|y_1|\le \n^{-\frac{1}{2}+\epsilon}$ and $|y_2|\ge \n^{-1+2\epsilon}$.
The above estimates imply that 
\beqq \begin{split}
	& \re \left( G(\gamma_1+\ii y_1, \wt \gamma_2+\ii y_2) - G(\gamma_1, \gamma_2) \right)
	\ge -c'  y_1^2 -\frac{1}4 \log(1+cy_2^2)
\end{split} \eeqq
and hence, the left-hand side of~\eqref{eq:y2cuff00} is bounded above by 
\beqq
	2 e^{-c' \n y_1^2} \int_{\n^{-1/2+2\epsilon}}^\infty e^{-\frac{\n}4 \log(1+cy_2^2)} \, \dd y_2 \,.
\eeqq	
The integral is uniformly bounded and $e^{-c' \n y_1^2} \le e^{-c'\n^{2\epsilon}}$. 
Hence we obtain the lemma.
\end{proof}

The above truncation is not enough. 
The next lemma show that we can truncate further to the interval $|y_2|\le  \n^{-\frac23+2\epsilon}$.
Here we use the fact that $(\gamma_1, \gamma_2)$ is the critical point. 
Note that this $y_2$-interval is smaller than the interval $|y_1|\le \n^{-\frac12+\epsilon}$.

\begin{lem} \label{lem:y_2 cutoff}
Uniformly for $|y_1| \leq \n^{-\frac{1}{2}+\epsilon}$, 
\beq \label{eq:cufo2} \begin{split}
	 \left( \int_{-\n^{-\frac12+\epsilon}}^{-\n^{-\frac23+2\epsilon}} + \int_{\n^{-\frac23+2\epsilon}}^{\n^{-\frac12+\epsilon}} \right)   \exp \left[ \n \left( G(\gamma_1+\ii y_1, \wt \gamma_2+\ii y_2) - G(\gamma_1, \gamma_2) \right) \right] \, \dd y_2
	= O(e^{-\n^{4\epsilon}})
\end{split} \eeq
\end{lem}

\begin{proof}
We start with~\eqref{eq:Gdiffincutf}. 
From the fact that $(\gamma_1,\gamma_2)$ is a critical point, we showed in~\eqref{critical low 00temp1} that 
\beq \label{eq:ratgammaB}
	\aln = \ben (\gamma_1-\gamma_2) .
\eeq
Inserting this into~\eqref{eq:Gdiffincutf} to remove $\aln$, and then expanding the terms involving $\ben$ in terms of powers of $y_1$, we find that 
\beq \begin{split} \label{G difference real}
	&\re \left( G(\gamma_1 + \ii y_1, \wt \gamma_2 + \ii y_2) - G(\gamma_1, \gamma_2) \right) \\
	&= - \frac{\ben (\gamma_1 + \gamma_2)}{2\gamma_1^2} y_1^2 - \frac{1}{4\n} \sum_{i=1}^{\n} \log \left[ \left( 1 -\frac{4y_1 y_2}{4\gamma_1 \gamma_2 - \mu_i} \right)^2 + \left( \frac{4\gamma_1 y_2}{4\gamma_1 \gamma_2 -\mu_i}  \right)^2 \right] + O(y_1^4).
\end{split} \eeq
From the rigidity, Condition \ref{cond:rigidity}, it is easy to check that 
\beq \label{eq:mu1N}
	\n^{-\frac{2}{3}} \ll \mu_1-\mu_{\n^{4\epsilon}} \ll \n^{-\frac{2}{3}+2\epsilon}. 
\eeq
The upper bound implies that $4\gamma_1 \gamma_2 - \mu_{\n^{4\epsilon}}  \ll \n^{-\frac{2}{3}+2\epsilon}$. 
Hence, for  $|y_2|\ge \n^{-\frac23+2\epsilon}$,
\beq \label{eq:cuffosusm} \begin{split}
	&\frac{1}{4\n} \sum_{i=1}^{\n^{4\epsilon}} \log \left[ \left ( 1 -\frac{4y_1 y_2}{4\gamma_1 \gamma_2 - \mu_i} \right)^2 + \left( \frac{4\gamma_1 y_2}{4\gamma_1 \gamma_2 -\mu_i}  \right)^2 \right]  \\
	&\geq \frac{1}{2\n} \sum_{i=1}^{\n^{4\epsilon}} \log \left( \frac{4 \gamma_1 |y_2|}{4\gamma_1 \gamma_2 - \mu_i} \right) \geq C\n^{-1+4\epsilon}.
\end{split} \eeq
The lower bound of~\eqref{eq:mu1N} implies that $4\gamma_1 \gamma_2 - \mu_{\n^{\epsilon}} \gg |y_1 y_2|$ for $|y_2|\le \n^{-\frac12+\epsilon}$ and $|y_1|\le \n^{-\frac12+\epsilon}$.
Hence, 
\beqq \begin{split}
	&\frac{1}{4\n} \sum_{i=\n^{4\epsilon}+1}^{\n} \log \left[ \left ( 1 -\frac{4y_1 y_2}{4\gamma_1 \gamma_2 - \mu_i} \right)^2 + \left( \frac{4\gamma_1 y_2}{4\gamma_1 \gamma_2 -\mu_i}  \right)^2 \right] 
	\geq \frac{1}{2\n} \sum_{i=\n^{4\epsilon}+1}^{\n} \log \left ( 1 -\frac{4y_1 y_2}{4\gamma_1 \gamma_2 - \mu_i} \right)  \\
	&
	\geq -\frac{C}{\n} \sum_{i=\n^{4\epsilon}+1}^{\n} \frac{|y_1 y_2|}{4\gamma_1 \gamma_2 - \mu_i} 
	\geq -C'|y_1y_2|
	\geq -C'\n^{-1+ 3\epsilon}.
\end{split} \eeqq
Note that the exponent $(-1+3\epsilon)$ is smaller than $(-1+4\epsilon)$ in~\eqref{eq:cuffosusm}.
Therefore, we obtain for $\n^{-\frac23+2\epsilon} \le |y_2|\le \n^{-\frac12+\epsilon}$ that
\beqq
	\re \left( G(\gamma_1 + \ii y_1, \wt \gamma_2 + \ii y_2) - G(\gamma_1, \gamma_2) \right) \leq -C \n^{-1+4\epsilon}.
\eeqq
This implies the lemma.
\end{proof}

\subsection{Decomposition of the double integral} \label{subsec:decompose}

We consider the part of the double integral~\eqref{double integral estimate 2} 
with $|y_1| \leq \n^{-\frac{1}{2} + \epsilon}$ and $|y_2| \leq \n^{-\frac{2}{3}+ 2\epsilon}$.
From~\eqref{eq:G diff 0101} and~\eqref{eq:ratgammaB},  using the Taylor series, 
\beq \begin{split} \label{G difference}
	&G(\gamma_1 + \ii y_1, \wt \gamma_2 + \ii y_2) - G(\gamma_1, \gamma_2) \\ 
	&= \ii \ben y_2 - \frac{\ben (\gamma_1 + \gamma_2)}{2\gamma_1^2} y_1^2 - \frac{1}{2\n} \sum_{i=1}^{\n} \log\left( 1 -\frac{4y_1 y_2}{4\gamma_1 \gamma_2 - \mu_i} + \frac{4\ii \gamma_1 y_2}{4\gamma_1 \gamma_2 - \mu_i} \right) + O(y_1^3).
\end{split} \eeq
Hence, 
\beq \begin{split} \label{G difference 200}
	&\exp \left[ \n\left( G(\gamma_1 + \ii y_1, \wt \gamma_2 + \ii y_2) - G(\gamma_1, \gamma_2) \right) \right] \\
	&\quad = \exp \left[ \ii \ben \n y_2 - \frac{\ben (\gamma_1 + \gamma_2)}{2\gamma_1^2} \n y_1^2 - \frac{1}{2} \sum_{i=1}^{\n} \log\left( 1 + \frac{4\ii \gamma_1 y_2}{4\gamma_1 \gamma_2 - \mu_i} \right) \right] \\
	&\qquad \qquad \times \exp \left[ -\frac{1}{2} \sum_{i=1}^{\n} \log\left( 1 - \frac{4y_1 y_2}{4\gamma_1 \gamma_2 - \mu_i + 4\ii \gamma_1 y_2} \right) + O(y_1^3) \right] .
\end{split} \eeq
Applying Lemma \ref{lem:logsummuigamma},
\beq
	\sum_{i=1}^{\n} \frac{1}{|4\gamma_1 \gamma_2 - \mu_i + 4\ii \gamma_1 y_2|^\ell} 
	\le  \sum_{i=1}^{\n} \frac{1}{|4\gamma_1 \gamma_2 - \mu_i|^\ell} 
	\le C_0^{\ell} \n^{\ell+\epsilon}
\eeq
for $\ell = 2, 3, \dots$, where $C_0$ is the constant in Lemma~\ref{lem:logsummuigamma}. For $\ell=1$, we use the bound
\beq
	\sum_{i=1}^{\n} \frac{1}{|4\gamma_1 \gamma_2 - \mu_i + 4\ii \gamma_1 y_2|} \leq \sum_{i=1}^{\n} \frac{1}{4\gamma_1 \gamma_2 - \mu_i} 
	= \frac{\n}{2\gamma_1} \left( \ben - \partial_2 G(\gamma_1, \gamma_2) \right) 
	= \frac{\n \ben}{2\gamma_1} = O(\n).
\eeq
This implies that, from the conditions on $y_1, y_2$, 
\beq
	\sum_{i=1}^{\n} \left( \frac{y_1y_2}{4\gamma_1 \gamma_2 - \mu_i + 4\ii \gamma_1 y_2} \right)^\ell
	= O(\n^{-\frac16 \ell+(3\ell+1)\epsilon})
\eeq
for $\ell = 2, 3, \dots$ and
\beq
	\sum_{i=1}^{\n} \frac{y_1y_2}{4\gamma_1 \gamma_2 - \mu_i + 4\ii \gamma_1 y_2}
	= O(\n^{-\frac16+3\epsilon}).
\eeq
Thus, expanding the last exponential function in~\eqref{G difference 200}, we obtain 
\beq \begin{split} \label{G difference 2}
	&\exp \left[ \n\left( G(\gamma_1 + \ii y_1, \wt \gamma_2 + \ii y_2) - G(\gamma_1, \gamma_2) \right) \right] \\
	& \quad = \exp \left[ \ii \ben \n y_2 - \frac{\ben (\gamma_1 + \gamma_2)}{2\gamma_1^2} \n y_1^2 - \frac{1}{2} \sum_{i=1}^{\n} \log\left( 1 + \frac{4\ii \gamma_1 y_2}{4\gamma_1 \gamma_2 - \mu_i} \right) \right] \\
	&\quad \qquad \qquad \times \left( 1 + 2y_1 y_2 \sum_{i=1}^{\n} \frac{1}{4\gamma_1 \gamma_2 - \mu_i + 4\ii \gamma_1 y_2} + O(\n^{-\frac{1}{3}+6\epsilon}) \right) .
\end{split} \eeq

We thus have
\beq
	\int_{-\n^{-\frac12+\epsilon}}^{\n^{-\frac12+\epsilon}}  \int_{-\n^{-\frac23+2\epsilon}}^{\n^{-\frac23+2\epsilon}}  \exp \left[ \n\left( G(\gamma_1 + \ii y_1, \wt \gamma_2 + \ii y_2) - G(\gamma_1, \gamma_2) \right) \right]  \, \dd y_2 \, \dd y_1 =: I_1 + I_2 + I_3,
\eeq
where $I_1$, $I_2$, and $I_3$ are given as follows: First, 
\beq \begin{split}
	I_1 &= \int_{-\n^{-\frac12+\epsilon}}^{\n^{-\frac12+\epsilon}}  \int_{-\n^{-\frac23+2\epsilon}}^{\n^{-\frac23+2\epsilon}}
	\exp \left[ \ii \ben \n y_2 - \frac{\ben (\gamma_1 + \gamma_2)}{2\gamma_1^2} \n y_1^2 - \frac{1}{2} \sum_{i=1}^{\n} \log\left( 1 + \frac{4\ii \gamma_1 y_2}{4\gamma_1 \gamma_2 - \mu_i} \right) \right] \dd y_2 \dd y_1. 
\end{split} \eeq
This is equal to the product of two single integrals $I_{11}$ and $I_{12}$. The $y_1$-integral is 
\beq \begin{split}
	 I_{11}:= \int_{-\n^{-\frac12+\epsilon}}^{\n^{-\frac12+\epsilon}} \exp \left[ - \frac{\ben (\gamma_1 + \gamma_2)}{2\gamma_1^2} \n y_1^2 \right] \dd y_1 .
\end{split} \eeq
This is real-valued and we have 
\beq
	\frac{C}{\sqrt{\n}}\le I_{11}\le \frac{C'}{\sqrt{\n}}.
\eeq
The $y_2$-integral is 
\beq \begin{split}
	I_{12} &:= \int_{-\n^{-\frac23+2\epsilon}}^{\n^{-\frac23+2\epsilon}}
	\exp \left[ \ii \ben \n y_2  - \frac{1}{2} \sum_{i=1}^{\n} \log\left( 1 + \frac{4\ii \gamma_1 y_2}{4\gamma_1 \gamma_2 - \mu_i} \right) \right] \dd y_2 . 
\end{split} \eeq
This is also real-valued since the imaginary part of the integrand is an odd function of $y_2$.
We have 
\beq \begin{split} \label{eq:I12ub}
	I_{12} \le  |I_{12}|&\le \int_{-\infty}^\infty 
	\exp \left[  - \frac{1}{4} \sum_{i=1}^{\n} \log\left( 1 + \left( \frac{4 \gamma_1 y_2}{4\gamma_1 \gamma_2 - \mu_i} \right)^2 \right) \right] \dd y_2  \\
	&\le  \int_{-\infty}^\infty 
	\exp \left[  - \frac{\n}{4}  \log\left( 1 + Cy_2^2 \right)  \right] \dd y_2  
	\le C.
\end{split} \eeq
On the other hand, we will show  the following lower bound in Subsection~\ref{subsec:I_12}: 
\beq \label{eq:I12esmc}
	 I_{12} \geq C\n^{-1-5\epsilon}. 
\eeq
Assuming this is true and using~\eqref{eq:I12ub}, we find that 
\beq \label{eq:Iallest}
	C\n^{-\frac32-5\epsilon} \le  I_{1}\le C \n^{-\frac12}.
\eeq

Second, 
\beqq \begin{split}
	I_2 &= \int_{-\n^{-\frac12+\epsilon}}^{\n^{-\frac12+\epsilon}}  \int_{-\n^{-\frac23+2\epsilon}}^{\n^{-\frac23+2\epsilon}}  \exp \left[ \ii \ben \n y_2 - \frac{\ben (\gamma_1 + \gamma_2)}{2\gamma_1^2} \n y_1^2 - \frac{1}{2} \sum_{i=1}^{\n} \log\left( 1 + \frac{4\ii \gamma_1 y_2}{4\gamma_1 \gamma_2 - \mu_i} \right) \right]  \\
	& \qquad \qquad \qquad \qquad \qquad \times \left( 2y_1y_2 \sum_{i=1}^{\n} \frac{1}{4\gamma_1 \gamma_2 - \mu_i + 4\ii \gamma_1 y_2} \right) \dd y_2 \dd y_1 \,.
\end{split} \eeqq
Since the integrand is an odd function of $y_1$, we find that $I_2 = 0$.

Finally, $I_3$ satisfies 
\beqq \begin{split}
	|I_3| &\leq C \n^{-\frac{1}{3}+6\epsilon}  \int_{-\n^{-\frac12+\epsilon}}^{\n^{-\frac12+\epsilon}}  \int_{-\n^{-\frac23+2\epsilon}}^{\n^{-\frac23+2\epsilon}}  \exp \left[  - \frac{\ben (\gamma_1 + \gamma_2)}{2\gamma_1^2} \n y_1^2 - \frac{1}{2} \re \sum_{i=1}^{\n}  \log\left( 1 + \frac{4\ii \gamma_1 y_2}{4\gamma_1 \gamma_2 - \mu_i} \right) \right] \dd y_2 \dd y_1 \\
	& \le C' \n^{-\frac{5}{6}+6\epsilon}   \int_{-\n^{-\frac23+2\epsilon}}^{\n^{-\frac23+2\epsilon}}  \exp \left[  - \frac{1}{2}  \log\left( \frac{4\gamma_1 y_2}{4\gamma_1 \gamma_2 - \mu_1} \right) \right] \dd y_2 \\
	&\leq C'' \n^{-\frac{5}{6}+6\epsilon} \int_{-\n^{-\frac{2}{3}+2\epsilon}}^{\n^{-\frac{2}{3}+2\epsilon}} \frac{1}{\sqrt{\n^{1-\epsilon} y_2}} \dd y_2  \leq C'' \n^{-\frac{5}{3} +8\epsilon}.
\end{split} \eeqq
Note that this upper bound is smaller than the lower bound of~\eqref{eq:Iallest} if $\epsilon < \frac{1}{78}$. 

Combining all estimates of Subsections \ref{subsec:truncate} and \ref{subsec:decompose}, we obtain
\beq \begin{split} \label{eq:lowtempmainines}
	C\n^{-\frac32-5\epsilon}\le - \int_{\gamma_1-\ii\infty}^{\gamma_1+\ii\infty} \int_{\gamma_2-\ii\infty}^{\gamma_2+\ii\infty} \exp \left[ \n \left( G(z_1, z_2) - G(\gamma_1, \gamma_2) \right) \right] \, \dd z_2 \, \dd z_1  \le C\n^{-\frac12} \, , 
\end{split} \eeq
thus prove Lemma~\ref{lem:lowt}, assuming that~\eqref{eq:I12esmc} is true. 

\subsection{Analysis of $I_{12}$} \label{subsec:I_12}

To complete the proof of Lemma~\ref{lem:lowt}, it remains to show the lower bound 
$I_{12} \geq C\n^{-1-5\epsilon}$ in~\eqref{eq:I12esmc}. 
We note by checking directly from the definition of $G$ that, 
\beq \begin{split}
	I_{12} &= \int_{-\n^{-\frac23+2\epsilon}}^{\n^{-\frac23+2\epsilon}}
	\exp \left[ \n\left( G(\gamma_1 , \gamma_2 + \ii y_2) - G(\gamma_1, \gamma_2) \right) \right] \dd y_2 . 
\end{split} \eeq
Define
\beq \label{K definition}
	K := -\ii \int_{\gamma_2-\ii\infty}^{\gamma_2+\ii\infty} \exp \left[ \n\left( G(\gamma_1, z) - G(\gamma_1, \gamma_2) \right) \right] \dd z.
\eeq
Then $I_{12}$ is the same integral as $K$ where the contour is restricted to the part $|z-\gamma_2|\le \n^{-\frac23+2\epsilon}$. 
Note that $K$ is real-valued since $G(\gamma_1, \overline{z_2})= \overline{G(\gamma_1, z_2)}$. 
The lower bound~\eqref{eq:I12esmc} follows if we show that 
\begin{enumerate}[(a)]
\item$|K- I_{12}|   \le e^{-C\n^{4\epsilon}}$, and 
\item $K\ge C\n^{-1-5\epsilon}$.
\end{enumerate}

Since $K$ is real-valued, 
\beqq
	|K-I_{12}|  
	\le \left(\int_{-\infty}^\infty- \int_{-\n^{-\frac23+2\epsilon}}^{\n^{-\frac23+2\epsilon}} \right)
	\exp \left[ \n\re \left( G(\gamma_1 , \gamma_2 + \ii y_2) - G(\gamma_1, \gamma_2) \right) \right] \dd y_2.
\eeqq
We have
\beqq \begin{split}
	\n\re \left( G(\gamma_1 , \gamma_2 + \ii y_2) - G(\gamma_1, \gamma_2) \right)
	&= -\frac14\sum_{i=1}^{\n} \log \left[ 1+ \left( \frac{4\gamma_1y_2}{4\gamma_1\gamma_2 -\mu_i} \right)^2 \right] .
\end{split}
\eeqq
Using~\eqref{eq:cuffosusm}, for $\n^{-\frac23+2\epsilon}\le |y_2|\le \n$, we have the estimate
\beqq \begin{split}
	\frac14\sum_{i=1}^{\n} \log \left[ 1+ \left( \frac{4\gamma_1y_2}{4\gamma_1\gamma_2 -\mu_i} \right)^2 \right] 
	&\ge  \frac14\sum_{i=1}^{\n^{4\epsilon}} \log \left[ 1+ \left( \frac{4\gamma_1y_2}{4\gamma_1\gamma_2 -\mu_i} \right]^2 \right)
	\ge C\n^{4\epsilon}.
\end{split}
\eeqq
For $|y_2|\ge \n$,
\beqq
	\frac14\sum_{i=1}^{\n} \log \left[ 1+ \left( \frac{4\gamma_1y_2}{4\gamma_1\gamma_2 -\mu_i} \right)^2 \right] 
	\ge \frac{\n}4 \log ( 1+ c y_2^2) \ge \frac{\n}4 \log (c|y_2|).
\eeqq
Hence
\beq
	|K-  I_{12}|  
	\le 2\n e^{ -C\n^{4\epsilon}} 
	+ 2\int_{\n}^\infty (c y_2)^{-\frac{\n}4} \dd y_2 \le e^{-C'\n^{4\epsilon}}. 
\eeq
We thus obtained property (a). 

We now prove the property (b), $K\ge C\n^{-1-5\epsilon}$.
We follow the proof of Lemma 6.3 in \cite{BaikLee} closely.
Observe that $\gamma_2$ is a critical point of the function $G(\gamma_1, z)$. 
Let $\Gamma$ be the curve of steepest-descent that passes through the point $\gamma_2$.
It satisfies $\im G(\gamma_1, z) = 0$.
It is straightforward to check from the formula of $G$ that 
\begin{enumerate}[(i)]
\item $\Gamma$ is symmetric about the real axis,
\item $\Gamma\cap \C^+$ is a $C^1$ curve,
\item $\Gamma$ lies in the half plane $\re z\le \gamma_2$, 
\item $\Gamma$ intersects with $\R$ only at $\gamma_2$,
\item the tangent line of $\Gamma$ at $\gamma_2$ is parallel to the imaginary axis,
\item the asymptote of $\Gamma$ is the negative real axis. 
\end{enumerate}

For example, the property (iii) can be checked by noting that for $z=x+\ii y$ with $x>\gamma_2$, 
\beqq
	F(y):= \im G(\gamma_1, x+\ii y)
	= \ben y - \frac1{2\n} \sum_{i=1}^{\n} \arctan\left( \frac{4\gamma_1y}{(4\gamma_1 x-\mu_i)^2+(4\gamma_1y)^2} \right)
\eeqq
has the global minimum at $y=0$ by computing its derivative. 

Since $\re (G(\gamma_1, z)-G(\gamma_1, \gamma_2))\le \ben \gamma_2 -\frac12 \log (R/2)$ for $|z|=R$ with $\re (z)\le \gamma_2$, 
we can deform the contour so that 
\beqq
	K = -\ii \int_{\Gamma} \exp \left[ \n\left( G(\gamma_1, z) - G(\gamma_1, \gamma_2) \right) \right] \dd z.
\eeqq
For $z \in \Gamma$, we let $x = \re z$ and $y = \im z$. Then, $\dd z = \dd x + \ii \dd y$ and
\beqq
	K= -\ii \int_{\Gamma} \exp \left[ \n\left( G(\gamma_1, z) - G(\gamma_1, \gamma_2) \right) \right] \dd x + \int_{\Gamma} \exp \left[ \n\left( G(\gamma_1, z) - G(\gamma_1, \gamma_2) \right) \right] \dd y.
\eeqq
Let $\Gamma^+ = \Gamma \cap \C^+$. By symmetry,
\beq
	K = 2 \int_{\Gamma^+} \exp \left[ \n\left( G(\gamma_1, z) - G(\gamma_1, \gamma_2) \right) \right] \dd y.
\eeq
In \cite{BaikLee}, the lower bound of $K$ was obtained by restricting the integral to a small ball of radius $\n^{-2}$. In the current work, however, we need to refine the argument further to prove \eqref{K definition}. We let $D_1$ be the disk of radius $\n^{-1-\epsilon}$ centered at $\gamma_2$, and similarly, $D_2$ be the disk of radius $\n^{-1-2\epsilon}$ centered at $\gamma_2$. The rest of the contour is controlled by the following lemma.

\begin{lem}[Lemma 6.4 of \cite{BaikLee}] \label{lem:decreasing}
Suppose that $f$ is a real-valued function defined on $\Gamma^+$ and $f(z)$ is decreasing along the curve $\Gamma^+$ as $z$ moves from the point $\gamma_2$ to the point $-\infty$. Then,
\beq
	\int_{\Gamma^+} e^{f(z)} \dd y \geq 0.
\eeq
\end{lem}

Since $G(\gamma_1, z)$ is analytic for $z_2$ in $D_1$, the series expansion 
\beq \label{G taylor}
	G(\gamma_1, z) - G(\gamma_1, \gamma_2) 
	= \sum_{j=2}^{\infty} \frac{1}{j!}  \partial_2^j G(\gamma_1, \gamma_2)(z-\gamma_2)^j
\eeq
converges for $z \in \Gamma^+ \cap D_1$. Set $X= \re (z-\gamma_2)$ and $Y=\im z = \im (z-\gamma_2)$. Comparing the imaginary parts of the both sides of \eqref{G taylor} by using Lemma \ref{lem:logsummuigamma}, we find that
\beq \label{XY curve}
	0 = \partial_2^2 G(\gamma_1, \gamma_2) XY + \frac{1}{2} \partial_2^3 G(\gamma_1, \gamma_2 X^2 Y - \frac{1}{6} \partial_2^3 G(\gamma_1, \gamma_2) Y^3 + \wt \Omega
\eeq
with
\beqq
	\wt \Omega = \sum_{j=4}^{\infty} \frac{1}{j!} \partial_2^j G(\gamma_1, \gamma_2) \im \left( (X+\ii Y)^j \right).
\eeqq
Note that $\im \left( (X+\ii Y)^j \right)$ is a homogeneous polynomial of $X$ and $Y$ with degree $j$. In the polynomial, every term contains both $X$ and $Y$, possibly except the term $Y^j$ when $j$ is odd. In any case,
\beqq
	\left| \frac{1}{j!} \partial_2^j G(\gamma_1, \gamma_2) \im \left( (X+\ii Y)^j \right) \right| \leq \frac{(8\gamma_1 C_0)^j}{j!} \n^{j-1+\epsilon} |z-\gamma_2|^{j-2} |XY| + \frac{(4 \gamma_1 C_0)^j}{j!} \n^{j-1+\epsilon} Y^j,
\eeqq
hence
\beq
	|\wt \Omega| \leq C\n^{1-\epsilon} (|XY| + Y^2).
\eeq

Define
\beqq
	\tau = - \frac{\partial_2^3 G(\gamma_1, \gamma_2)}{  \partial_2^2 G(\gamma_1, \gamma_2)} > 0.
\eeqq
From Lemma~\ref{lem:logsummuigamma} (by putting $\epsilon/4$ instead of $\epsilon$), we find that 
$\partial_2^2 G(\gamma_1, \gamma_2) \geq C\n^{1-\frac{\epsilon}{2}}$ and 
\beq
	\tau \ll \n^{1-\epsilon}. 
\eeq
Thus, dividing both sides of \eqref{XY curve} by $\partial_2^2 G(\gamma_1, \gamma_2) Y$, we obtain that
\beqq
	X(1+o(1)) + \frac{\tau}{6} Y^2 (1+ o(1)) = 0,
\eeqq
and thus, 
\beq \label{eq:XY}
	X = -\frac{\tau}{6} Y^2 (1+o(1)).
\eeq
We also see that $\Gamma^+ \cap D_1$ is a graph, $\dd y = \dd Y$ is positive on $\Gamma^+ \cap D_1$, and $\Gamma^+$ intersects $\partial D_1$ at exactly one point.

Let $\zeta_1$ (resp. $\zeta_2$) be the point where $\Gamma^+$ and $\partial D_1$ (resp. $\partial D_2$) intersect. Then,
\beq \begin{split} \label{zeta_1 estimate}
	G(\gamma_1, \gamma_2) - G(\gamma_1, \zeta_1) &\geq \frac{1}{4} \partial_2^2 G(\gamma_1, \gamma_2) |\zeta_1 - \gamma_2|^2 - \sum_{j=3}^{\infty} \frac{1}{j!} C_0^j \n^{j-1+\frac{\epsilon}{4}} |\zeta_1 - \gamma_2|^j \\
	&\geq C\n^{-1-\frac{5\epsilon}{2}}
\end{split} \eeq
and
\beq \begin{split} \label{zeta_2 estimate}
	G(\gamma_1, \gamma_2) - G(\gamma_1, \zeta_2) &\leq \frac{1}{2} \partial_2^2 G(\gamma_1, \gamma_2) |\zeta_2 - \gamma_2|^2 + \sum_{j=3}^{\infty} \frac{1}{j!} C_0^j \n^{j-1+\frac{\epsilon}{4}} |\zeta_2 - \gamma_2|^j \\
	&\leq C\n^{-1-\frac{7\epsilon}{2}}.
\end{split} \eeq
We introduce the function
\beqq
	f(z) = \begin{cases}
	\n G(\gamma_1, \zeta_1) - G(\gamma_1, \gamma_2) & \text{if } z \in \Gamma^+ \cap D_1, \\
	\n G(\gamma_1, z) - G(\gamma_1, \gamma_2) & \text{if } z \in \Gamma^+ \cap D_1^c.
	\end{cases}
\eeqq
It is obvious that $f(z)$ is a decreasing function of $z$ along the curve $\Gamma^+$ as $z$ moves from $\gamma_2$ to $-\infty$. Thus, applying Lemma \ref{lem:decreasing} to the function $f$,
\beq \label{bound on D_1}
\int_{\Gamma^+ \cap D_1} \exp \left[ \n\left( G(\gamma_1, \zeta_1) - G(\gamma_1, \gamma_2) \right) \right] \dd y + \int_{\Gamma^+ \cap D_1^c} \exp \left[ \n \left( G(\gamma_1, z) - G(\gamma_1, \gamma_2) \right) \right] \dd y \geq 0.
\eeq
Since $\dd y$ is positive on $\Gamma^+ \cap D_1$,
\beq \begin{split} \label{K split}
	\frac{K}{2} &= \int_{\Gamma^+} \exp \left[ \n\left( G(\gamma_1, z) - G(\gamma_1, \gamma_2) \right) \right] \dd y \\
	&= \int_{\Gamma^+ \cap D_2} \exp \left[ \n\left( G(\gamma_1, z) - G(\gamma_1, \gamma_2) \right) \right] \dd y + \int_{\Gamma^+ \cap D_2^c} \exp \left[ \n\left( G(\gamma_1, z) - G(\gamma_1, \gamma_2) \right) \right] \dd y \\
	&\geq \int_{\Gamma^+ \cap D_2} \exp \left[ \n\left( G(\gamma_1, \zeta_2) - G(\gamma_1, \gamma_2) \right) \right] \dd y + \int_{\Gamma^+ \cap D_2^c} \exp \left[ \n\left( G(\gamma_1, z) - G(\gamma_1, \gamma_2) \right) \right] \dd y.
\end{split} \eeq
Subtracting \eqref{bound on D_1} from \eqref{K split}, we find that
\beq \begin{split}
	\frac{K}{2} &\geq \int_{\Gamma^+ \cap D_2} \big( \exp \left[ \n\left( G(\gamma_1, \zeta_2) - G(\gamma_1, \gamma_2) \right) \right] - \exp \left[ \n\left( G(\gamma_1, \zeta_1) - G(\gamma_1, \gamma_2) \right) \right] \big) \dd y \\
	&\qquad + \int_{\Gamma^+ \cap D_1 \cap D_2^c} \big( \exp \left[ \n\left( G(\gamma_1, z) - G(\gamma_1, \gamma_2) \right) \right] - \exp \left[ \n\left( G(\gamma_1, \zeta_1) - G(\gamma_1, \gamma_2) \right) \right] \big) \dd y \\
	&\geq \int_{\Gamma^+ \cap D_2} \big( \exp \left[ \n\left( G(\gamma_1, \zeta_2) - G(\gamma_1, \gamma_2) \right) \right] - \exp \left[ \n\left( G(\gamma_1, \zeta_1) - G(\gamma_1, \gamma_2) \right) \right] \big) \dd y.
\end{split} \eeq
From \eqref{eq:XY}, we find that $\im \zeta_2 \gg \re \zeta_2$, hence $\im \zeta_2 \geq C\n^{-1-2\epsilon}$. Since $\int_{\Gamma^+ \cap D_2} \dd y = \im \zeta_2 \geq C\n^{-1-2\epsilon}$, we find from the estimates \eqref{zeta_1 estimate} and \eqref{zeta_2 estimate} that
\beq
	\frac{K}{2} \geq C\n^{-1-2\epsilon} \left[ \exp (-C\n^{-\frac{7\epsilon}{2}}) - \exp (-C\n^{-\frac{5\epsilon}{2}}) \right] \geq C\n^{-1-5\epsilon}.
\eeq
This proves the desired lower bound of $K$. Thus, Lemma~\ref{lem:lowt} is proved.

\subsection{Double integral in the low temperature regime} \label{sec:lowtempresult}

The following is the main result for the random double integral in the low temperature regime.

\begin{prop}[Random double integral for low temperature] \label{prop:lowresult}
Assume Conditions~\ref{cond:aanndb}~\ref{cond:empirical} and~\ref{cond:rigidity}.
Suppose that $\be$ in Condition~\ref{cond:aanndb} satisfies $\be>\be_c$ where $\be_c$ is defined in \eqref{eq:B_c}. 
Then, for every $\epsilon>0$, 
\beq \begin{split}
	\frac1{\n} \log \q= 
	& \ew + (\mu_1-\red) \lw+ O(\n^{-1+\epsilon}) 
\end{split} \eeq
where
\beq \begin{split}
	\ew &= \sqrt{\al^2+\red \be^2} - \al \log \left(  \frac{\al + \sqrt{\al^2+\red \be^2}}{2\be} \right)  -\frac12 \h(\red), \\
	\lw &= \frac{\be^2}{2(\al + \sqrt{\al^2+\red \be^2})} -\frac12 \h'(\red).
\end{split} \eeq
\end{prop}

\begin{proof}
From Lemma~\ref{lem:lowt},  
\beqq \begin{split}
	\frac1{\n} \log \q= G(\gamma_1, \gamma_2)  + O\left( \frac{\log \n}{\n} \right) .
\end{split} \eeqq
We have
\beq \begin{split} \label{eq:theothertwoterms}
	G(\gamma_1, \gamma_2)  
	=\sqrt{\aln^2+\gamma \ben^2} - \frac1{2\n} \sum_{i=1}^{\n} \log(\gamma-\mu_i)
	- \aln \log \left(  \frac{\aln + \sqrt{\aln^2+\gamma \ben^2}}{2\ben} \right) .
\end{split} \eeq
From Lemma~\ref{lem:logsummuigamma},
\beqq \begin{split}
	\frac1{\n} \sum_{i=1}^{\n} \log(\gamma-\mu_i)
	= \h(\red)+(\mu_1-\red)\h'(\red) + O(\n^{-1+\epsilon}).
\end{split} \eeqq
On the other hand, for the other two terms in~\eqref{eq:theothertwoterms}, we first replace 
$\aln$ and $\ben$ by $\al$ and $\be$ and introduce an error term $O(\n^{-1-\delta})$. 
We then replace $\gamma$ by $\mu_1$ and introduce an error term $O(\n^{-1+\epsilon})$ due to Lemma~\ref{lem:lowcp}. 
Writing $\mu_1= \red + (\mu_1-\red)$ and using the Taylor expansion up to the first order using $(\mu_1-\red)^2= O(\n^{-4/3+2\epsilon})$, we find that 
\beqq \begin{split}
	&\sqrt{\aln^2+\gamma \ben^2} - \aln \log \left(  \frac{\aln + \sqrt{\aln^2+\gamma \ben^2}}{2\ben} \right) \\
	&=\sqrt{\al^2+\red \be^2} - \al \log \left(  \frac{\al + \sqrt{\al^2+\red \be^2}}{2\be} \right) 
	+ \frac{\be^2}{2(\al + \sqrt{\al^2+\red \be^2})} (\mu_1-\red) + O(\n^{-4/3+2\epsilon}).
\end{split} \eeqq
This completes the proof. 
\end{proof}

\section{Proof of Theorem~\ref{thm:eigval}} \label{sec:proofeigval}

Recall that we assume $N_1\ge N_2$ and $N=N_1+N_2$. We take the limit as $N, N_1, N_2\to \infty$ satisfying 
\beq
	\frac{N_1}{N} = \rf + O(N^{-1-\delta}), \qquad \frac{N_2}{N}= \rs+ O(N^{-1-\delta})
\eeq
for $\rf\ge \rs>0$ and $\rf+\rs=1$. 
From~\eqref{eq:partition_Q}, we have
\beq \label{eq:ZQRa}
	Z_{N_1,N_2}(\beta) 
	= \qq \bigg( N_2, \frac{N_1-N_2}{2N_2}, \frac{N_1}{\sqrt{N_2N}} \beta \bigg) R(N_1, N_2, \beta)
\eeq
where $\qq$ is defined with the eigenvalues $\mu_i$, $1\le i\le N_2$, of the random matrix $\frac1{N_1} J^TJ$. 
Here $J$ is an $N_1\times N_2$ matrix with independent and identically distributed entries of mean $0$ and variance $1$ satisfying the assumptions in Subsection~\ref{sec:defn}. 
The constant 
\beq
	R(N_1, N_2, \beta)
	=  \frac{1}{|S^{N_1-1}||S^{N_2-1}|} 2^{N_2} \left( \frac{\pi^2 N}{N_1^2N_2 \beta^2} \right)^{(N-4)/4} .
\eeq
We use the results of the previous two sections on $\q$. 
To do that, we need to check Conditions~\ref{cond:aanndb},~\ref{cond:empirical} and~\ref{cond:rigidity}:
\begin{itemize}
\item Conditions~\ref{cond:aanndb} holds with  
\beq \begin{split} \label{eq:abpar}
	&\n = N_2, \qquad \al=\frac{\rf-\rs}{2\rs}, \qquad \be = \frac{\rf}{\sqrt{\rs}} \beta. 
\end{split} \eeq
where $\aln= \frac{N_1-N_2}{2N_2}$ and $\ben= \frac{N_1}{\sqrt{N_2N}} \beta$.
\item Condition~\ref{cond:empirical} follows from the well-known Marchenko-Pastur law~\cite{Marcenko-Pastur63} in random matrix theory.
The limiting empirical measure is given by~\eqref{eq:MPformula},  
\beq \label{eq:redled}
	\dd \wh \mu(x)= \dd\mu_{\MP}(x) := \frac{2 \sqrt{(\red-x)(x-\led)}}{\pi (\sqrt{\red}-\sqrt{\led})^2 x} \mathds{1}_{(\led, \red)}(x) \dd x
\eeq
where
\beq
	\led= \frac{(\sqrt{\rf}-\sqrt{\rs})^2}{\rf}, \qquad \red= \frac{(\sqrt{\rf}+ \sqrt{\rs})^2}{\rf} .
\eeq
\item Condition~\ref{cond:rigidity} holds with high probability for the eigenvalues. This is proved recently by Pillai and Yin in~\cite{PY}  for $\rf > \rs$. For the case $\rf=\rs$, Condition~\ref{cond:rigidity} with high probability follows from Corollary 1.3 of~\cite{AjankiErdosKruger2014} and the fact that $\dd \mu_{\MP}(x) = \dd \mu_{\SC}(\sqrt x)$ where $\mu_{\SC}$ is the Wigner semicircle distribution. (See also Equation (1.12) of \cite{AjankiErdosKruger2014}.)
Similar rigidity results were proved for other classes of random matrices starting with the Wigner matrices~\cite{EYY} and also various random matrix models including invariant ensembles~\cite{BEY} and sparse random matrices~\cite{EKYY1}. Rigidity estimates are obtained from the local laws such as local semicircle law or local Marchenko-Pastur law, and they are also crucial a priori estimates for the proof of bulk and edge universality of random matrices.
\end{itemize}

\bigskip 
Before we deduce the limit of $\qq$, we first state the asymptotics of $R(N_1, N_2, \beta)$.

\begin{lem} \label{lem:Rcom}
We have 
\beq \begin{split}
	 \frac1N\log \left( R(N_1, N_2, \beta) \right) = 
	 &-\frac 12 + \rs\log 2 - \frac12 \log(2\rf\sqrt{\rs} \beta) + \frac{\rf}{2} \log \rf + \frac{\rs}{2} \log \rs \\
	&  \quad 
	 + \frac{\log N}{N} - \frac1{N}  \log \left( \frac{\pi}{\rf \sqrt{\rf\rs} \beta^2} \right) 
	 + O(N^{-2}) .
\end{split} \eeq
\end{lem}

\begin{proof}
The area of unit sphere satisfies 
\beq
	\log (|S^{n-1}|)
	= \log \left( \frac{2\pi^{n/2}}{\Gamma(n/2)} \right)
	=  \frac{n}{2} \log \left(\frac{2\pi e}{n}\right) + \frac{1}{2} \log \left( \frac{n}{\pi} \right) + O(n^{-1}).
\eeq
Hence, 
\beqq \begin{split}
	\frac1N\log (|S^{N_1-1}||S^{N_2-1}|) 
	&= \frac12 \log \left(\frac{2\pi e}{N}\right) 
	- \frac{\rf}{2} \log \rf - \frac{\rs}{2} \log \rs 
	+ \frac1{N} \log \left( \frac{\sqrt{\rf\rs} N}{\pi} \right) + O(N^{-2}).
\end{split} \eeqq
On the other hand, 
\beqq \begin{split}
	\frac1N\log \left( 2^{N_2} \left( \frac{\pi^2 N}{N_1N_2^2 \beta^2} \right)^{(N-4)/4} \right)
	= \rs\log 2 + \left( \frac12 -\frac{2}{N} \right) \log \left( \frac{\pi}{\rf \sqrt{\rs} \beta N} \right).
\end{split} \eeqq
We thus obtain the lemma.
\end{proof}

\subsection{Transforms of the Marchenko-Pastur distribution and critical temperature}

We need the following formulas of the log transform and the Stieltjes transform of the Marchenko-Pastur distribution. 

\begin{lem} \label{lem:MPform}
Set 
\beq \label{eq:Rfistd}
 	R(z):= \sqrt{(z-\led)(z-\red)}. 
\eeq
We have 
\beqq
\begin{split}
	H_{\MP}(z) &:=\int_{\R} \log(z-x) \dd\mu_{\MP}(x) \\
	&= \frac{2}{(\sqrt{\red}-\sqrt{\led})^2} \bigg[ z- \sqrt{\red\red} \log z- R(z) +(\red+\led)\log(\sqrt{z-\led}+\sqrt{z-\red})  \\
	&\qquad\qquad  + \sqrt{\red\led} \log \left( \frac{ \sqrt{\red(z-\led)}-\sqrt{\led(z-\red)} }{ \sqrt{\red(z-\led)}+\sqrt{\led(z-\red)} } \right)
	 \bigg] 
\end{split}
\eeqq
for $z\notin (-\infty,  \red)$
and
\beq \label{eq:sMP}
	s_{\MP}(z) :=H_{\MP}'(z) =\int_{\R} \frac1{z-x} \dd\mu_{\MP}(x) 
	= \frac{2 (z-R(z) -  \sqrt{\red\red}) }{(\sqrt{\red}-\sqrt{\led})^2 z}
\eeq
for $z\notin (\led,  \red)$.
\end{lem}

\begin{proof}
The computation of the Stieltjes transform~\eqref{eq:sMP} is a standard exercise in complex analysis. 
The log transform $H_{\MP}(z)$ is obtained from~\eqref{eq:sMP} by taking anti-derivative. 
\end{proof}

In terms of $\rf$ and $\rs$,  
\beq \label{eq:Rhere}
 	R(z) = \frac{\sqrt{(\rf z- 1)^2-4\rf\rs}}{\rf},
\eeq
\beq
\begin{split}
	H_{\MP}(z)
	 &= \frac{1}{2\rs} \bigg[ \rf z - 1 - \rf R(z) - (\rf-\rs) \log z+ \log \left( \frac{\rf z- 1+\rf R(z)}{2 \rf} \right)  \\
	&\qquad\qquad + (\rf-\rs) \log \left( \frac{  \rf z- (\rf-\rs)\rf  R(z)- (\rf-\rs)^2 }{ 2\rf \rs  z } \right)  	 \bigg] ,
\end{split}
\eeq
and
\beq \label{eq:sMPhp}
	s_{\MP}(z) =H_{\MP}'(z) 
	= \frac{\rf z- \rf R(z) -  ( \rf-\rs ) }{2 \rs  z} .
\eeq

\bigskip

Let us consider the critical inverse temperature $\beta_c$. 
From~\eqref{eq:abpar}, $\beta_c= \frac{\sqrt{\rs}}{\rf} B_c$, and from Definition~\eqref{defn:cv}, 
$\be_c= \sqrt{\red (s_{\MP}(\red))^2+2\al s_{\MP}(\red)}$. 
From the above explicit formulas, we see that 
\beq
	s_{\MP}(\red) = H_{\MP}'(\red)= \frac{\rf}{\sqrt{\rs}(\sqrt{\rf}+\sqrt{\rs})}. 
\eeq
Hence we find that the critical inverse temperature is 
\beq
	\beta_c= (\rf\rs)^{-\frac14}. 
\eeq

\subsection{High temperature case}

For $\beta<\beta_c$, we evaluate the terms in Proposition~\ref{prop:integral high} explicitly. 
We first find $z_c$ solving the equation~\eqref{eq:cpeq}. 
From~\eqref{eq:abpar} and~\eqref{eq:sMPhp}, this equation is 
\beq \label{eq:H_MP'}
	\frac{\rf z_c- \rf R(z_c)- (\rf-\rs)}{2\rs z_c} = \frac{2\rf^2\beta^2}{\rf-\rs+ \sqrt{W(z_c)}}
\eeq
where $R(z)$ is given by~\eqref{eq:Rhere} and 
\beq
	W(z):= (\rf-\rs)^2+4\rf^2\rs \beta^2 z .
\eeq
We claim that the solution of this equation in $z_c\in (\red, \infty)$ is given by 
\beq \label{eq:z_c solution}
	z_c = \frac{(1+\rf\beta^2)(1+\rs\beta^2)}{\rf\beta^2} = \frac{1+\beta^2+\rf\rs\beta^4}{\rf\beta^2}.
\eeq
Indeed, with this $z_c$, we find that 
\beq \label{eq:Rzcf}
	\sqrt{W(z_c)}= 1+ 2 \rf\rs\beta^2 , \qquad R(z_c)= \frac{1-\rf\rs\beta^4}{\rf\beta^2}
\eeq
where we used the condition that $\beta<(\rf\rs)^{-1/4}$.  
Using this, we find that both sides of~\eqref{eq:H_MP'} are $\frac{\rf\beta^2}{1+\rs\beta^2}$.
This verifies~\eqref{eq:z_c solution}.

Now from Proposition~\ref{prop:integral high} and recalling that $\n=N_2$, we find that for any $\epsilon > 0$,
\beq \begin{split}
	&\frac1{N}\log \qq = \at - \frac{1}{2N} \left[  \sum_{i=1}^{N_2} \log\left ( z_c - \mu_i\right) - N_2 H_{\MP}(z_c) \right] 
	- \frac{\log (\rs N)}{N} + \frac{1}{2N} \log \left( \frac{\pi^2}{\dw}  \right) + O(N^{-1-\epsilon}) 
\end{split} \eeq
with high probability where
\beq \begin{split}
	\at=\rs\aw &= \frac{\sqrt{W(z_c)}}{2} - \frac{\rf-\rs}{2} \log \left(  \frac{\rf-\rs+ \sqrt{W(z_c)}}{4\rf\sqrt{\rs}\beta}\right)  -\frac{\rs}2 H_{\MP}(z_c)
\end{split} \eeq
and
\beq \begin{split}
	\dw &=  - \frac{4(\rf-\rs)}{\rs} H_{\MP}''(z_c) -8z_c H_{\MP}'(z_c) H_{MP}''(z_c) -4 (H_{\MP}'(z_c))^2 .
\end{split} \eeq

It is direct to check that
\beq \begin{split}
	H_{\MP}(z_c) 
	= \rf \beta^2 + \frac{\rf-\rs}{\rs} \log \left( \frac{1}{1+ \rs\beta^2} \right) - \log \left( \rf\beta^2\right)
\end{split} \eeq
and
\beq
	H_{\MP}'(z_c) =\frac{\rf\beta^2}{1+\rs\beta^2}, 
	\qquad
	H_{\MP}''(z_c) = - \frac{\rf^2\beta^4}{(1+\rs\beta^2)^2(1-\rf\rs\beta^4)}.
\eeq
Thus, we obtain  
\beq
	\at = \frac12(1+\rf\rs\beta^2) - \frac{\rf-\rs}2 \log \left( \frac{1}{2\sqrt{\rs}\beta} \right)  
	+ \frac{\rs}2 \log \left( \rf\beta^2\right) 
\eeq
and 
\beq
	\dw = \frac{4\rf^3\beta^4}{\rs(1-\rf\rs\beta^4)}.
\eeq
Recalling~\eqref{eq:ZQRa} and using Lemma~\ref{lem:Rcom}, we find that for $\beta<\beta_c$,
\beq \begin{split}
	\frac1{N} \log Z_{N_1,N_2}(\beta) 
	 = \frac{\rf\rs\beta^2}2 - \frac{1}{2N} &\left[  \sum_{i=1}^{N_2} \log\left ( z_c - \mu_i\right) - N_2 H_{\MP}(z_c) \right]  \\
	+ \frac1{N} &\left[  \frac12 \log \left(1-\rf\rs\beta^4\right) -\log 2 \right]
	 + O(N^{-1-\epsilon}) 
\end{split} \eeq
with high probability. This proves the part (i) of Theorem~\ref{thm:eigval}.


\subsection{Low temperature case}

For $\beta<\beta_c$, Proposition~\ref{prop:lowresult} implies that (recall that $\n=N_2$ and $\frac{N_2}{N}\to \rs$)
\beq \begin{split}
	\frac1{N} \log \qq(\n, \al, \be) = 
	& \et + (\mu_1-\red) \lt+ O(N^{-1+\epsilon}).
\end{split} \eeq
where
\beq \begin{split}
	\et = \rs\ew &= \frac{\sqrt{W(\red)}}{2} - \frac{\rf-\rs}{2} \log \left(  \frac{\rf-\rs+ \sqrt{W(\red)}}{4\rf\sqrt{\rs}\beta}\right)  -\frac{\rs}2 H_{\MP}(\red),  \\
	\lt = \rs\lw  &= \frac{\rf^2\rs\beta^2}{\rf-\rs+ \sqrt{W(\red)}} -\frac{\rs}2 H_{\MP}'(\red)
\end{split} \eeq
with 
\beq
	W(\red)
	= (\sqrt{\rf}+\sqrt{\rs})^2((\sqrt{\rf}-\sqrt{\rs})^2+4\rf\rs\beta^2) .
\eeq
From the explicit formulas, 
\beq \begin{split}
	H_{\MP}(\red) 
	&= \sqrt{\frac{\rf}{\rs}} - \frac{\rf-\rs}{\rs} \log \left( \frac{\sqrt{\rf}+\sqrt{\rs}}{\sqrt{\rf}} \right) 
	+  \log \sqrt{\frac{\rs}{\rf}} 
\end{split} \eeq
and
\beq
	H_{\MP}'(\red)= \frac{\rf}{\sqrt{\rs}(\sqrt{\rf}+\sqrt{\rs})} .
\eeq
Setting
\beq
	\varl=\varl(\rf, \rs):= (\sqrt{\rf}-\sqrt{\rs})^2+4 \rf\rs\beta^2 ,
\eeq	
we have
\beq \begin{split}
	\et	= \frac{(\sqrt{\rf}+\sqrt{\rs})\sqrt{\varl}-\sqrt{\rf\rs}}{2} - \frac{\rf-\rs}2 \log \left( \frac{\sqrt{\varl}+\sqrt{\rf}-\sqrt{\rs}}{4\sqrt{\rf\rs}\beta} \right)  - \frac{\rs}4 \log \frac{\rs}{\rf}
\end{split} \eeq
and 
\beq \begin{split}
	\lt = \frac{\rf\sqrt{\rf\rs}}{2(\sqrt{\rf}+\sqrt{\rs})} \left( \frac{2 \sqrt{\rf\rs}\beta^2 }{\sqrt{S}+ \sqrt{\rf}-\sqrt{\rs}} -\frac1{\sqrt{\rf}} \right)=\frac{ \rf(\sqrt{\varl} - \sqrt{\rf}-\sqrt{\rs}) }{ 4(\sqrt{\rf}+\sqrt{\rs})} . 
\end{split} \eeq
Using $\left( \frac{\sqrt{A^2+B^2}+A}{B}\right)^2= \frac{\sqrt{A^2+B^2}+A}{\sqrt{A^2+B^2}-A}$, we write the $\log$ term in $\et$ as 
\beq \begin{split}
	 \log \left( \frac{\sqrt{\varl}+\sqrt{\rf}-\sqrt{\rs}}{4\sqrt{\rf\rs}\beta} \right)  
	 = \frac12 \log \left( \frac{\sqrt{\varl}+\sqrt{\rf}-\sqrt{\rs}}{\sqrt{\varl}-\sqrt{\rf}+\sqrt{\rs}} \right)- \log 2 .
\end{split} \eeq
Hence, recalling~\eqref{eq:ZQRa} and using Lemma~\ref{lem:Rcom}, we find that for $\beta>\beta_c$, 
\beqq \begin{split} 
	&\frac1{N} \log Z_{N_1,N_2}(\beta) \\
	&=  \frac{(\sqrt{\rf}+\sqrt{\rs})\sqrt{\varl}-\sqrt{\rf\rs}-1}{2} - \frac{\rf-\rs}4 \log \left( \frac{\sqrt{\varl}+\sqrt{\rf}-\sqrt{\rs}}{\sqrt{\varl}-\sqrt{\rf}+\sqrt{\rs}} \right)  - \frac{\rs}4 \log \rf - \frac{\rf}{4} \log \rs - \frac12 \log \beta\\
	&\quad +   
	\left(\mu_1- \frac{(\sqrt{\rf}+\sqrt{\rs})^2}{\rf} \right) \frac{\rf ( \sqrt{\varl} - \sqrt{\rf}-\sqrt{\rs})}{4(\sqrt{\rf}+\sqrt{\rs})} + O(N^{-1+\epsilon})
\end{split} \eeqq
with high probability. This completes the proof of the part (ii) of Theorem~\ref{thm:eigval}.

\section{Proof of Theorem~\ref{thm:main}} \label{sec:ProofFlucFromRMT}

We derive Theorem~\ref{thm:main} from Theorem~\ref{thm:eigval}.

\subsection{High temperature}

When $\beta<\beta_c$, from Theorem~\ref{thm:eigval}, we need the behavior of the sum of $\log (z_c-\mu_i)$. 
It follows form the following result in random matrix theory. 

\begin{prop}[Linear statistics of the eigenvalues] \label{cond:linear}
For any function $\varphi:\R\to \R$ that is analytic in an open neighborhood of the support of $\mu_{\MP}$, the random variable
\beq
	\sum_{i=1}^{N_2} \varphi(\mu_i) - N_2 \int \varphi(x) \dd \mu_{\MP}(x) 
\eeq
converges in distribution as $N_2\to\infty$ to a Gaussian random variable with mean $M(\varphi)$ and variance $V(\varphi)$ given as follows: setting 
\beq
	\Phi(x)=\varphi \left ( \frac{\red-\led}4 x+ \frac{\red+\led}2 \right),
\eeq	
we have 
\beq \label{eq:MVfml}
	M(\varphi) = M_{\GOE}(\Phi)- (W_4-3) \tau_2(\Phi), 
	\qquad V(\varphi)=V_{\GOE}(\Phi) + (W_4-3) \tau_1(\Phi)^2
\eeq
where
\beq \begin{split}
	M_{\GOE}(\Phi) &= \frac{\Phi(-2)+\Phi(2)}{4} -\frac{\tau_0(\Phi)}{2}, \\
	V_{\GOE}(\Phi) &= \frac1{2\pi^2} \int_{\led}^{\red} \int_{\led}^{\red} \left( \frac{\Phi(x_1)- \Phi(x_2)}{x_1-x_2} \right)^2 \frac{4-x_1x_2}{\sqrt{4-x_1^2}\sqrt{4-x_2^2}} \, \dd x_1 \dd x_2.
\end{split} \eeq
and for $\ell=0,1,2,\cdots$, 
\beq
	\tau_{\ell}(\Phi)= \frac1{2\pi} \int_{-\pi}^\pi \Phi(2\cos\theta) \cos(\ell \theta) \dd \theta .
\eeq
\end{prop}

This result was first obtained in \cite{BS2004} (equation (5.13)). It was also obtained in \cite{LP} (Theorem 4.5) and \cite{BWZ2010} (Theorem 1.1). 
In the above formulas, GOE stands for Gaussian orthogonal ensemble, another classical random matrix ensemble.

To complete the proof of the part (i) of Theorem~\ref{thm:main}, we need to evaluate $M(\varphi)$ and $V(\varphi)$ 
when $\varphi(x)=\log (z_c-x)$. In this case, 
\beq
	\Phi(y)= \log \left( \frac{\red-\led}{4} \right) + \log (\tilde{z}-y), \qquad 
	\tilde{z}= \frac{4}{\red-\led} \left( z_c- \frac{\red+\led}2 \right).
\eeq
Since $\tau_0(1)=1$, we find that 
\beq
	\tau_\ell(\Phi)= \delta_{\ell=0}  \log \left( \frac{\red-\led}{4} \right) + \tau_{\ell}(\Phi_0), 
	\qquad \Phi_0(x):=\log(\tilde{z}-x).
\eeq
We computed in Appendix A of \cite{BaikLee} that for $\Phi_0(x)=\log(\tilde{z}-x)$, 
\beqq \begin{split}
	\tau_0(\Phi_0)= \log(\tilde{z}+\sqrt{\tilde{z}^2-4}) -\log 2, 
	\quad 
	\tau_1(\Phi_0)= \frac{\sqrt{\tilde{z}^2-4}}2- \frac{\tilde{z}}{2}, 
	\quad \tau_2(\Phi_0) = \frac{\tilde{z} \sqrt{\tilde{z}^2-4}}{4} - \frac{\tilde{z}^2}{4} + \frac{1}{2} .
\end{split} \eeqq
Note that
\beq
	\sqrt{\tilde{z}^2-4} = \frac{4}{\red-\led} \sqrt{(z_c-\led)(z_c-\red)} = \frac{4}{\red-\led} R(z_c),
\eeq
which can be checked from the definition~\eqref{eq:Rfistd}. 
Hence,
\beq \begin{split}
	 M_{\GOE}(\Phi) &= \frac{\log  R(z_c) }{2} -\frac{1}{2} \left[ \log\left( z_c- \frac{\red+\led}2+ R(z_c) \right) -\log 2 \right] , \\
	\tau_1(\Phi_0) &=  \frac{2}{\red-\led} \left[ R(z_c)  - \left( z_c- \frac{\red+\led}2 \right) \right] ,\\
	\tau_2(\Phi_0) &= \frac{4\left( z_c- \frac{\red+\led}2 \right)}{(\red-\led)^2} \left[ R(z_c) - \left( z_c- \frac{\red+\led}2 \right) \right] + \frac{1}{2}  .
\end{split} \eeq
On the other hand, from Lemma A.1 of \cite{BaikLee},
\beq
	V_{\GOE}(\Phi)= 2 \log \left( \frac{z_c- \frac{\red+\led}{2}+ R(z_c) }{2 R(z_c) } \right) .
\eeq
Inserting~\eqref{eq:redled} for $\red$ and $\red$, and recalling~\eqref{eq:Rzcf} for $R(z_c)$, we find that   
\beq \begin{split}
	&\tau_1(\Phi_0)=  -\sqrt{\rf\rs}\beta^2,  \qquad
	\tau_2(\Phi_0) = - \frac{\rf\rs\beta^4}{2} \\
	&M_{\GOE}(\Phi) =\frac{1}{2}  \log \left(  1-\rf\rs\beta^4 \right) , 
	\qquad V_{\GOE}(\Phi)= - 2 \log \left( 1-\rf\rs\beta^4 \right).
\end{split} \eeq
Hence, we have 
\beq \begin{split}  \label{V_phi}
	M(\varphi) &= \frac{1}{2}  \log \left(  1-\rf\rs\beta^4 \right) + (W_4-3) \frac{\rf\rs\beta^4}{2}, \\
	V(\varphi) &= - 2 \log \left( 1-\rf\rs\beta^4 \right) + (W_4-3) \rf\rs\beta^4 .
\end{split} \eeq
for $\varphi(x)=\log (z_c-x)$.
 
Part (i) of Theorem~\ref{thm:eigval}, Prosition~\ref{cond:linear} and~\eqref{V_phi} prove the part (i) of Theorem~\ref{thm:main} when $\rf\ge \rs$. 
The case when $\rf<\rs$ follows from the symmetry, noting that $F(\beta)$, $\mea$ and $\van$ are all symmetric in $\rf$ and $\rs$.

\subsection{Low temperature}

When $\beta<\beta_c$, we need the behavior of the top eigenvalue $\mu_1$. 
The following result is well-known in random matrix theory. 
See, for example, \cite{Johnstone01, Soshnikov2002} and also Corollary 1.2 of \cite{PY}. 
	
\begin{prop}[Tracy--Widom limit of the largest eigenvalue] \label{cond:edge}
We have
\beq
	\frac{N_1}{\sqrt{N_1}+\sqrt{N_2}} \left( \frac{1}{\sqrt{N_1}} + \frac{1}{\sqrt{N_2}} \right)^{-\frac{1}{3}}  (\mu_1 - \red)
	\Rightarrow \TW
\eeq
in distribution. 
\end{prop}

In terms of $\rf, \rs$, 
\beq
	\frac{N_1}{\sqrt{N_1}+\sqrt{N_2}} \left( \frac{1}{\sqrt{N_1}} + \frac{1}{\sqrt{N_2}} \right)^{-\frac{1}{3}}
	= \frac{\rf (\rf\rs)^{\frac{1}{6}}}{(\sqrt{\rf}+\sqrt{\rs})^{\frac{4}{3}}} N^{2/3} (1+ O(N^{-1-\delta})).
\eeq
Then, combining with the constant in~\eqref{eq:Ffllo}, 
\beq \begin{split}
	\frac{(\sqrt{\rf}+\sqrt{\rs})^{\frac{4}{3}}}{\rf (\rf\rs)^{\frac{1}{6}}} \frac{\rf ( \sqrt{\varl} - \sqrt{\rf}-\sqrt{\rs})}{4(\sqrt{\rf}+\sqrt{\rs})}	= \frac{ (\sqrt{\rf}+\sqrt{\rs})^{1/3} ( \sqrt{\varl} - \sqrt{\rf}-\sqrt{\rs}) }{ 4(\rf\rs)^{1/6} }  .
\end{split} \eeq
This proves the part (ii) of Theorem~\ref{thm:main} when $\rf\ge \rs$. 
The case when $\rf<\rs$ again follows from the symmetry noting that $F(\beta)$ and $\varial$ are symmetric in $\rf$ and $\rs$.

\section{Non-identically distributed disorders} \label{sec:non-identical}

In this section, we briefly discuss the case where the disorders are non-identically distributed. 
Let $\Sigma$ be a positive-definite matrix of size $N_1\times N_1$. 
Let $J$ be the matrix of i.i.d. entries as before.
Consider the new Hamiltonian 
\beq \label{eq:Hmt}
	H(\bss, \tt) = \frac1{\sqrt {N}} \sum_{i=1}^{N_1} \sum_{j=1}^{N_2} (\Sigma^{1/2}J)_{ij}\sigma_i \tau_j . 
\eeq 
The new disorder parameters are $(\Sigma^{1/2}J)_{ij}=\sum_k \Sigma^{1/2}_{ik} J_{kj}$. 
In particular, when $\Sigma^{1/2}$ is a diagonal matrix, the variances of the disorder parameters depend on the index $i$ but not on $j$. 

The associated random matrix is 
\beq
	S = \frac{1}{N_1} J^T \Sigma J.
\eeq
In statistics, $S$ is known as a sample covariance matrix with general population covariance matrix $\Sigma$.
Recall the ingredients of the analysis we have done in this paper. The double integral formula in Lemma~\ref{lem:representation}, which was the starting point of our analysis, holds for any sample covariance matrices: we set $M=\frac1{\sqrt{N_1}}\Sigma^{1/2}J$ in the proof. The proofs of asymptotic formulas for $\q$ in Proposition~\ref{prop:integral high} and Proposition~\ref{prop:lowresult} require the regularity of measure in Condition~\ref{cond:empirical} and the rigidity of the eigenvalues in Condition~\ref{cond:rigidity}. Finally, the fluctuation of the free energy was obtained by applying the linear statistics of the eigenvalues in Proposition~\ref{cond:linear} and the Tracy--Widom limit of the largest eigenvalue in Proposition~\ref{cond:edge}.

The limiting empirical spectral distribution (ESD) of $S$ is well studied in random matrix theory. Under a very general assumption, it was proved in \cite{SilversteinChoi} that the limiting ESD of $S$ is regular. The typical assumption is as follows: Let $\sigma_1 \geq \sigma_2 \geq \cdots \geq \sigma_{N_1}$ be the eigenvalues of $\Sigma$ and denote by $\widehat\sigma$ the empirical spectral distribution of $\Sigma$. Define $\xi_+$ as the unique solution in $(0, \sigma_1^{-1})$ to the equation
\beq
	\int \left( \frac{t \xi_+}{1-t \xi_+} \right)^2 \dd \widehat\sigma(t) = \frac{\rs^2}{\rf^2}.
\eeq
If
\beq
	\limsup \sigma_1 < \infty, \quad \liminf \sigma_{N_1} > 0, \quad \limsup \sigma_1 \xi_+ < 1,
\eeq
then the limiting ESD of $S$ is regular, and hence satisfy Condition~\ref{cond:empirical}. (See Theorem 3.1 of \cite{BPZ2015}.) 
This assumption basically means that the largest eigenvalue $\sigma_1$ of $\Sigma$ is not too far away from the rest of the eigenvalues.
Such an assumption was also used in \cite{El_Karoui06, LS_sample, KY_sample}. Under the same assumption, it was proved by Knowles and Yin \cite{KY_sample} that Condition~\ref{cond:rigidity} holds with high probability.

Proposition~\ref{cond:linear}, the central limit theorem for the linear statistics of the eigenvalues was first proved by Bai and Silverstein \cite{BS2004} when $W_4$, the fourth moment of $J_{ij}$, matches that of the standard Gaussian, which is $3$. It was later extended to a general case with any finite fourth moment by Najim and Yao \cite{NajimYao}.

Proposition~\ref{cond:edge} was first proved for the complex case with Gaussian disorder $J_{ij}$ by El Karoui \cite{El_Karoui06}. Bao, Pan, and Zhou \cite{BPZ2015} proved the edge universality for the model, which asserts that the rescaled distribution of the largest eigenvalue does not depend on the distribution of $J_{ij}$. The result in \cite{BPZ2015} also holds for the real case, but the Tracy--Widom limit for the case was not proved. Proposition~\ref{cond:edge} was later proved in \cite{LS_sample, KY_sample}.

Therefore, the main theorems, Theorem~\ref{thm:main} and Theorem~\ref{thm:eigval}, also hold for the Hamiltonian with non-identically distributed disorders under a general assumption on the spectrum of the matrix $\Sigma$, with suitable changes on the constants $\mu$, $\sigma^2$, $\varial$, and $\varl$.


\def\cydot{\leavevmode\raise.4ex\hbox{.}}

\end{document}